\documentclass[12pt,a4paper]{article}

\usepackage{pdflscape}
\usepackage{amsmath}
\usepackage{amssymb}
\usepackage{latexsym}
\usepackage{srcltx}
\usepackage{graphics}
\usepackage{color}
\usepackage{epsfig}
\usepackage{color}

\usepackage{subcaption}

\usepackage{anyfontsize}
\usepackage{tabularx}

\usepackage{graphicx}
\usepackage{multirow}
\usepackage{array}

\usepackage{amssymb, amsmath, amsfonts, mathtools}
\usepackage{esint}
\usepackage{graphicx}
\usepackage{breqn}
\usepackage{float}
\usepackage{caption}
\usepackage{subcaption}
\usepackage{hhline}
\usepackage{multirow}
\usepackage{rotating}
\usepackage{fixltx2e}

\newcommand{\nid}{\noindent}

\newcommand{\X}{X}

\newcommand{\cF}{{\mathcal{F}}}

\newcommand{\cO}{{\mathcal{O}}}

\newcommand{\cB}{{\mathcal{B}}}

\newcommand{\cH}{{\mathcal{H}}}
\newcommand{\cS}{{\mathcal{S}}}

\newcommand{\cQ}{{\mathcal{Q}}}

\newcommand{\cM}{{\mathcal{M}}}

\newcommand{\bx}{{\boldsymbol x}}

\newcommand{\ba}{{\boldsymbol a}}

\newcommand{\be}{{\boldsymbol e}}

\newcommand{\bv}{{\boldsymbol v}}
\newcommand{\bu}{{\boldsymbol u}}

\newcommand{\sgn}{{\rm sgn}}

\newtheorem{theorem}{Theorem}
\newtheorem{prop}{Proposition}
\newtheorem{lemma}{Lemma}
\newtheorem{cor}{Corollary}

\newtheorem{defi}{Definition}

\date{}
\author{Aleksandar Cvetkovi\' c\thanks{Gran Sasso Science Institute (GSSI), L'Aquila, Italy}}

\title{Stabilising the Metzler matrices with applications to dynamical systems}

\begin{document}

\maketitle

\begin{abstract}

Metzler matrices play a crucial role in  positive linear dynamical systems. Finding the closest stable Metzler matrix to an unstable one (and vice versa) is an important issue with many applications. The stability considered here is in the sense of Hurwitz, and the distance between matrices is measured in  $l_\infty,\  l_1$,  and in the max norms. We provide either explicit solutions or efficient algorithms for obtaining the closest (un)stable matrix. The procedure for finding the closest stable Metzler matrix is based on the recently introduced selective greedy spectral method for optimizing the Perron eigenvalue. Originally intended for non-negative matrices, here is generalized to Metzler matrices. The efficiency of the new algorithms is demonstrated in examples and by numerical experiments in the dimension of up to 2000. Applications to dynamical systems, linear switching systems, and sign-matrices are considered. 

\bigskip

\smallskip

\begin{flushright}
\noindent  \textbf{AMS 2010} {\em subject
classification: 15A42, 15B35, 15B48, 90C26, 93C30}
\end{flushright}

\end{abstract}
\bigskip

\begin{center}
\large{\textbf{1. Introduction}}
\end{center}

\medskip

Metzler matrices and the problem of finding their closest Hurwitz (un)stable counterparts arise in linear dynamical systems, differential equation analysis, electrodynamics, population dynamics, economics, etc. \cite{FaRi,Lue,Mit,And, Log}. We refer to the problems of finding the closest (un)stable matrix as the \textit{(de)stabilization problems}. By the \textit{closest matrix} we mean a matrix $X$ which minimizes the distance from the starting matrix $A$, in some matrix norm. In this paper we are concerned with $l_\infty,\ l_1$ and $\max$ norms. In other norms, such as Euclidean or Frobenius norms, the stabilization problem is hard and  as a rule does not allow finding global solutions \cite{And, GuPr}.\\

The analogous problems for the non-negative matrices and Schur stability were considered in \cite{AkGa,GiSh, GuPr,NePr}. In \cite{NePr} an explicit solutions were given for the problems of destabilization. For the stabilization, an efficient algorithms were developed. Of special interest is the algorithm for finding the closest stable non-negative matrix in $l_\infty$ and $l_1$ norms, based on strikingly efficient \textit{greedy spectral simplex method}. This method was devised in the same paper, and it is used for minimizing the spectral radius on the product families of (sparse) non-negative matrices. The greedy method can also be of independent interest, apart from the stabilization problem. Its modification, \textit{selective greedy method}, was recently developed in \cite{PrCv}. It is significantly easier to implement and  has a  quadratic convergence.\\

In order to stabilize the Metzler matrix, we need to deal with the problem of optimizing spectral abscissa on families of Metzler families. In general, this problem is notoriously hard, since the objective function is neither convex nor concave, nor is Lipschitz. Because of this, there might be many points of local extrema which are, moreover, hard to identify. However, for matrix families with special structure, i.e. \textit{product (or uncertainty) structure}, it is possible to find an efficient algorithm for optimizing the spectral abscissa. In this paper we present such algorithm, based on selective greedy method.\\

Furthermore, we discuss the application of stabilization problem for \textit{finding the closest stable Metzler sign-matrix}. Metzler sign-matrices and sign-stability are interesting concepts with many applications \cite{Bri,QuRu,Rob,May,Cla}, and here we present the algorithm for solving the sign-stabilization problem. Building upon this, we propose a procedure for stabilizing a positive linear switching system (LSS) \cite{Lib,GuSho} having a certain structure.\\

The paper is structured as follows:\\

\nid In Section 2 we deal with the Hurwitz (de)stabilization in the $\max$-norm, on the set of Metzler matrices. An explicit solution for the destabilization is given, while for the stabilization an iterative procedure is provided.\\

\nid Section 3 generalizes both the greedy and selective greedy methods for optimizing Perron eigenvalue on Metzler product families. They remain as effective as when used on non-negative matrices. This is demonstrated by the numerical results: even for the dimension of 2000, the computational time does not exceed 3 minutes (on a standard laptop).\\

\nid Section 4 is concerned with Hurwitz (de)stabilization in the $l_\infty$ norm\footnote{\ See the beginning of the next subsection.}. Again, an explicit solution for the closest unstable Metzler matrix is obtained. As for the stabilization, we provide an iterative algorithm. With slight modification, this algorithm can be applied for Schur stabilization of non-negative matrices. This modification significantly improves the computing time for Schur stabilization from \cite{NePr}.\\

\nid The applications for Hurwitz stabilization (in $l_\infty$ norm) are presented in Section 5.

\bigskip

\begin{center}
\textbf{1.1.  Some preliminaries}
\end{center}

\medskip

 Throughout this paper we assume, unless otherwise is stated, that all our matrices are $d\times d$ square. The matrix norms considered are $l_\infty,\ l_1$ and max norms, given by \cite{HoJo}:
\begin{equation*}
\begin{array}{lcc}
\|X\|_\max  & = & \max\limits_{1\leqslant i,j\leqslant d}\ |x_{ij}|;
\vspace*{0.05cm}\\
\|X\|_1 & =  & \max\limits_{1\leqslant j\leqslant d}\ \sum_{i=1}^d |x_{ij}|;
\vspace*{0.05cm}\\
\|X\|_\infty & =  & \max\limits_{1\leqslant i\leqslant d}\ \sum_{j=1}^d |x_{ij}|.
\end{array}
\end{equation*}
\nid As it can be easily seen $\|X\|_\infty = \|X^T\|_1$. Consequently, all the results for the $l_\infty$ norm apply for $l_1$ norm as well: we just need to take the transpose of the matrix. Having this in mind, we shall only develop results for $l_\infty$ and max norms\footnote{\ Norms without a subscript will be regarded as an $l_\infty$ norms (i.e $\|\cdot\|\ \equiv\|\cdot\|_\infty$).}. We now define the following sets of matrices:

\smallskip

\begin{defi}
A matrix $A$ is said to be {\rm Metzler} if all its off-diagonal elements are non-negative. We denote the set of all Metzler matrices by $\cM$. 
\end{defi} 

\begin{defi}
A Metzler matrix is {\rm non-negative} if all of its elements are non-negative. A non-negative matrix is {\rm strictly positive} if it contains no zero entry. If $A$ is a non-negative matrix we write $A\geqslant 0$, and if $A$ is strictly positive we write $A > 0$.
\end{defi}

\nid For two matrices $A$ and $B$ we write $A\geqslant B$ if $A - B\geqslant 0$, and $A\leqslant B$ if $B - A\geqslant 0$. Analogously, we define strict relations $A > B$ and $A < B$. 

\begin{defi}
A Metzler matrix is {\rm strict} if it contains a strictly negative entry.
\end{defi}

\begin{defi}
A Metzler matrix is {\rm full} if all of its off-diagonal elements are strictly positive.
\end{defi}

\smallskip

\nid An apparent property of a strict Metzler matrix, which we will amply exploit, is that it can be translated to a non-negative matrix: i.e. if $A$ is a strict Metzler matrix, then there exists $h > 0$ such that $A + hI$ is non-negative. Here (and throughout the paper) $I$ denotes the identity matrix.  

\bigskip

The spectral radius and spectral abscissa of a Metzler matrix are fundamental concepts in the development of our work, which we define below:

\smallskip

\begin{defi}
{\rm Spectral radius} of a Metzler matrix $A$, denoted by $\rho(A)$, is the largest modulus of its eigenvalues, that is:
$$\rho(A) = \max\{\ |\lambda|\ |\ \lambda {\rm\ is\ an\ eigenvalue\ of\ } A\in\cM\ \}.$$
\end{defi}

\begin{defi}
{\rm Spectral abscissa} of a Metzler matrix $A$, denoted by $\eta(A)$, is largest real part of all of its eigenvalues, that is:
$$\eta(A) = \max\{\ {\rm Re}(\lambda)\ |\ \lambda {\rm\ is\ an\ eigenvalue\ of\ } A\in\cM\ \}.$$
\end{defi}

\begin{defi}
The {\rm leading eigenvalue} (or {\rm Perron eigenvalue}) of a Metzler matrix is the eigenvalue with the largest real part, and its corresponding eigenvector is the {\rm leading (Perron) eigenvector}.
\end{defi}

\smallskip

\nid From the Perron-Frobenius theory we know that a non-negative matrix has its spectral radius as a leading eigenvalue \cite{Gan}. Since this implies that the leading eigenvalue is real and non-negative, we have:

\begin{prop}\label{etaro}
If $A$ is a non-negative matrix, then
$$\rho(A) = \eta(A).$$
{\hfill $\Box$}
\end{prop}

\smallskip

\nid Specially, a strictly positive matrix has a simple leading eigenvalue, and its corresponding leading eigenvector is unique and strictly positive. For the non-negative matrix, this might not always be the case: its leading eigenvalue might not be simple, its leading eigenvector can contain zero entries and not necessarily be unique.\\ 

Given the fact that any strict Metzler matrix $A$ can be translated to a non-negative one $A+hI$, for some $h>0$, and that both $A$ and $A+hI$ have the same sets of eigenvectors, from Proposition \ref{etaro} we obtain:

\begin{lemma}\label{trans}
If $A$ is strict Metzler, then
$$\rho(A+hI) = \eta(A) + h,$$
for any $h>0$ such that $A+hI$ is non-negative.
{\hfill $\Box$}
\end{lemma}

It is well known that a spectral radius is monotone on the set of non-negative matrices \cite{NePr}. Now, we prove that the same holds for the spectral abscissa on the set of Metzler matrices, that is:

\begin{lemma}\label{monotonija}
Let $A$ and $B$ be two Metzler matrices such that $B\geqslant A$. Then $\eta(B)\geqslant\eta(A)$.
\end{lemma} 

\noindent {\tt Proof.} Let $h>0$ be such that both $A + hI$ and $B + hI$ are non-negative. Since $B$ is entrywise bigger than $A$ we have $B + hI\geqslant A + hI$. By the monotonicity of spectral radius $\rho(B + hI)\geqslant\rho(A + hI)$ holds. Using Lemma \ref{trans}, we obtain $\eta(B)\geqslant\eta(A)$. {\hfill $\Box$}\\

Even though the results of the Perron-Frobenius theory regarding the spectral radius do not hold for the strict Metlzer matrices, the following Theorem \ref{trans.tm} reveals that, in this case, we should just change our focus to the spectral abscissa. 

\begin{theorem}{\rm\cite{BeFa}}\label{trans.tm}
The full strict Metzler matrix has its spectral abscissa as a (simple) leading eigenvalue, its corresponding leading eigenvector is unique and strictly positive. A strict Metzler matrix has also its spectral abscissa as a leading eigenvalue, although it may not be simple, its corresponding leading eigenvector is non-negative and may not necessarily be unique.
{\hfill $\Box$}
\end{theorem}

\medskip

\smallskip

\begin{defi}
A Metzler matrix $A$ is {\rm strongly Schur stable} if $\rho(A) < 1$, and {\rm weakly Schur stable} if $\rho(A)\leqslant 1$. Otherwise, if $\rho(A) > 1$ we say it is {\rm strongly Schur unstable}, and if $\rho(A)\geqslant 1$ we say it is {\rm weakly Schur unstable}.
\end{defi}

\begin{defi}
A Metzler matrix $A$ is {\rm strongly Hurwitz stable} if $\eta(A) < 0$, and {\rm weakly Hurwitz stable} if $\eta(A)\leqslant 0$. Otherwise, if $\eta(A) > 0$ we say it is {\rm strongly Hurwitz unstable}, and if $\eta(A)\geqslant 0$ we say it is {\rm weakly Hurwitz unstable}.
\end{defi}

\smallskip

\nid We denote by $\cH$ the set of all weakly Hurwitz stable Metzler matrices, and by $\cH_s$ the set of all strongly Hurwitz stable Metzler matrices.\\

\medskip

% \nid \textbf{Remark.} Even though the concept of Schur stability is valid for any % Metzler matrix, we will be concerned about Schur stability only in the case of non-% negative matrices. This is because spectral radius does not play a role of the Perron %eigenvalue on the set of (strict) Metzler matrices.\\

\nid \textbf{Remark.} When searching for the closest (un)stable Metzler matrix to the matrix $A$, the starting matrix need not necessarily be Metzler. It is easy to check that both the real matrix $A$ and matrix $A'$, with entries

\begin{equation*}
a'_{ij} = 
\left\{
\begin{array}{cc}
a_{ij}, & a_{ij}\geqslant 0\ {\rm or}\ i = j \\
0, & {\rm otherwise}
\end{array}
\right.
\end{equation*}\\

\nid have the same solution. Therefore, without the loss of generality, we shall always assume that our starting matrix is Metzler.\\

\bigskip

\begin{center}\label{sec2}
\large{\textbf{2. Problems of (de)stabilization in the max-norm}}
\end{center}

\medskip

\begin{center}
\textbf{2.1. Closest Hurwitz unstable Metzler matrix in the max-norm}
\end{center}

\medskip

For a given Hurwitz stable Metzler matrix $A$ we consider the problem of finding its closest (weakly) Hurwitz unstable matrix $X$, with respect to the max-norm. In other words, for a matrix $A\in\cH_s$ find a matrix $X$ that satisfies:
\begin{equation*}
\left\{
\begin{array}{l}
\|X - A\|_{\max}\ \rightarrow\min\\
\eta(X) = 0,\ X\in\cM.
\end{array}
\right.
\end{equation*}
The following set of results provides the answer to the posed problem. We start from\\

\begin{lemma}{\rm\cite{NePr}}\label{lem1}
A non-negative matrix $A$ is Schur stable if and only if the matrix $(I - A)^{-1}$ is well defined and non-negative.
{\hfill $\Box$}
\end{lemma}

\noindent and its analogue for Metzler matrices: 

\begin{lemma}{\rm\cite{Bri}}\label{lem2} 
A Metzler matrix $A$ is Hurwitz stable if and only if it is invertible and the matrix $-A^{-1}$ is non-negative.
{\hfill $\Box$}
\end{lemma}

\begin{lemma}\label{lem3} Let $A\in\cH_s$ and $H\geqslant 0$, $H\neq 0$. Then

\begin{equation}\label{eqn1}
\eta(A + \alpha H) < 0,\qquad 0\leqslant\alpha < \frac{1}{\rho(-A^{-1}H)}
\end{equation}
and $\eta\Big(A + \frac{H}{\rho(-A^{-1}H)}\Big) = 0$.
\end{lemma}

\noindent {\tt Proof.} Define matrix $B = -A^{-1}H$. Since $A$ is Hurwitz stable, from Lemma \ref{lem2} we have $B\geqslant 0$. Define
$$W(\alpha) = -A - \alpha H = -A - \alpha(-AB) = -A(I - \alpha B).\vspace*{0.1cm}$$
From Lemma \ref{lem1} and Lemma \ref{lem2}, $W(\alpha)$ is invertible for every $\alpha\in\Big[0,\frac{1}{\rho(B)}\Big)$. Moreover, for\vspace*{0.1cm} every $\alpha$ from the given interval $W^{-1}(\alpha)$ is non-negative. Since $-W(\alpha) = A + \alpha H$ is Metzler, by Lemma \ref{lem2}, it is Hurwitz stable, so we have (\ref{eqn1}).\\

For the proof of the second part, we assume that matrix $B$ is strictly positive. Then its leading eigenvector $\bv$ is also strictly positive, and we have $-W\Big(\frac{1}{\rho(B)}\Big)\bv = 0$. Hence, by the continuity of spectral abscissa on the set of Metzler matrices, we have   $\eta\Big(A + \frac{H}{\rho(B)}\Big) = 0$.
{\hfill $\Box$}\\

\medskip

\nid The following result is a direct consequence of Lemma \ref{lem3}:

\begin{cor}\label{this.cor}
Let $A\in\cH_s$ and $H\geqslant 0,\ H\neq 0$. Then all $X\in\cM$ that satisfy

$$X < A + \frac{H}{\rho(-A^{-1}H)}$$
are stable.
{\hfill $\Box$}
\end{cor}

\noindent Applying Corollary \ref{this.cor} to the special case of the non-negative matrix $H = E$, where $E$ is matrix of all ones, we arrive to the explicit formula for Hurwitz destabilization in the max-norm.

\begin{theorem}
Let $A$ be a Hurwitz stable Metzler matrix and $\be$ vector of all ones. Then all Metzler matrices $X$ that satisfy
$$X < A - \frac{E}{\langle A^{-1}\be,\be\rangle}$$
are Hurwitz stable. Moreover, matrix $A - \frac{E}{\langle A^{-1}\be,\be\rangle}$ is the closest Hurwitz unstable Metzler\vspace*{0.1cm} matrix with the distance to the starting matrix $A$ equal to
$$\tau_* = \frac{1}{\sum\limits_{i,j}|a'_{ij}|},$$
where $a'_{ij}$ are entries of the matrix $A^{-1}$. 
\end{theorem}

\noindent {\tt Proof.} The proof stems directly from the fact that
$$\rho(-A^{-1}E) = -\langle A^{-1}\be,\be\rangle.$$
{\hfill $\Box$}

\bigskip

\begin{center}
\textbf{2.2. Closest Hurwitz stable Metzler matrix in the max-norm}
\end{center}

\medskip

We now deal with a problem of Hurwitz stabilization in the max-norm. First, consider the following spectral abscissa minimization problem for a given matrix $A\in\cM$ and parameter $\tau\geqslant 0$:
\begin{equation}\label{uneq.1}
\min\limits_{X\in\cM}\{\eta(X)\ |\ \|X-A\|_{\max}\leqslant\tau\}.
\end{equation}

\begin{lemma}\label{opt.stab.max}
The optimal solution of the problem {\rm(\ref{uneq.1})} is a matrix $A(\tau)\in\cM$ with the following entries:
\begin{equation*}
a_{ij}(\tau) = \left\{
\begin{array}{l}
a_{ii} - \tau,\ i=j\\
\max\{0,a_{ij}-\tau\},\ i\neq j.
\end{array}
\right.
\end{equation*}
\end{lemma}

\nid {\tt Proof.} Clearly, a Metzler matrix $A(\tau)$ is feasible for (\ref{uneq.1}). Moreover, for any other feasible solution $X$, we have $X\geqslant A(\tau)$. From the monotonicity of spectral abscissa we have that $A(\tau)$ is truly the optimal solution to our problem.
{\hfill $\Box$}\\

Now, for a given Hurwitz unstable Metzler matrix $A$, consider the following minimization problem:
\begin{equation}\label{uneq.2}
\min\limits_{X\in\cM}\|X-A\|_{\max}.
\end{equation}

\nid Denote by $\tau_*$ its optimal value. 

\begin{lemma}\label{unlem}
$\tau_*$ is a unique root of the equation $\eta(A(\tau)) = 0$.
\end{lemma}
\nid {\tt Proof.} The statement follows directly from  Lemma \ref{opt.stab.max} and the fact the function $\eta(A(\tau))$ is monotonically decreasing in $\tau$.
{\hfill $\Box$}\\

\nid The matrix $A(\tau_*)$ is the closest Hurwitz stable Metzler matrix to a given matrix $A\in\cM$ with $\eta(A)>0$, and $\tau_*$ is the minimal distance. We now present the algorithm for computing the value $\tau_*$ and solving (\ref{uneq.2}).\\

\bigskip

\nid {\tt Algorithm 1: Finding the closest Hurwitz stable Metzler matrix in max-norm}

\bigskip

\nid {\tt Step 1.} Sort all positive entries of matrix $A$ in an increasing order. Check if the highest entry is on the main diagonal. If so, take $\tau_* = \max_i a_{ii}$ and finish the procedure. Else, continue to the\\

\nid {\tt Step 2.} Using the monotonicity of the spectral abscissa, find by bisection in the entry number the value $\tau_1$, which is the largest between all positive entries $a_{ij}$ and zero, such that $\eta(A(a_{ij}))>0$. Find value $\tau_2$, which is the smallest entry of $A$ with $\eta(A(a_{ij}))<0$.\\

\nid {\tt Step 3.} Form the matrix
$$H = \frac{A(\tau_1) - A(\tau_2)}{\tau_2 - \tau_1}.$$\\

\nid {\tt Step 4.} Compute
$$\tau_* = \tau_2 - \frac{1}{\rho(-A^{-1}(\tau_2)H)}.$$
{\hfill $\diamond$}

\bigskip

\begin{theorem}
Algorithm 1 computes an optimal solution to the problem {\rm (\ref{uneq.2})}.
\end{theorem}

\nid {\tt Proof.} 
First, assume that the highest entry is on the main diagonal. Then, the matrix $A(\max_i a_{ii})$ will be non-positive diagonal matrix such that $\eta(A(\max_i a_{ii})) = 0$. Since $\tau_*$, by Lemma \ref{unlem}, is the unique root of the equation $\eta(A(\tau)) = 0$, we can put $\tau_* = \max_i a_{ii}$.\\  

Now, assume that the highest entry is off-diagonal. Then, the matrix $A(\max_{ij}a_{ij})$ is negative diagonal. Therefore, we have
$$\eta\Big(A\Big(\max_{ij}a_{ij}\Big)\Big)<0.$$
On the other hand, $\eta(A(0))=\eta(A)>0$. So, it is possible to find two values $\tau_1<\tau_2$ from the set 
$$\{0\}\cup\{a_{ij}\ |\ a_{ij}>0\}$$
such that $\eta(A(\tau_1))>0$ and $\eta(A(\tau_2))<0$. $A(\tau)$ is linear in $\tau$, and for every $\tau\in[\tau_1,\tau_2]$ we have 
$$A(\tau) = A(\tau_2) + \frac{\tau_2 - \tau}{\tau_2 - \tau_1}(A(\tau_1) - A(\tau_2)) = A(\tau_2) + (\tau_2 - \tau)H.$$
Applying Lemma \ref{lem3} for Metzler matrix $A(\tau_2)$ and non-negative matrix $H$ concludes the proof.
{\hfill $\Box$}

\bigskip

\begin{center}\label{sec3}
\large{\textbf{3. Optimizing the spectral abscissa on product families of Metzler matrices}}
\end{center}

\medskip

Stabilization in $l_\infty$ norm is much more delicate than in the max-norm. In order to properly address it, we first need to find a way to optimize spectral abscissa on Metzler product families\footnote{\ For the reminder of the text, each time we mention \textit{Metzler product families} or \textit{positive product families}, or similar, it will be in the context of \textit{product families of Metzler matrices}, \textit{product families of positive matrices}, etc.}. Optimizing spectral abscissa on such families can also be of independent interest, since product families occur quite often in applications \cite{BloNe,NePr2,Pro,CaLipp,Log,Koz}.\\

By $\cM_i$ we denote the sets of vectors
$$\cM_i = \{\ \ba\in\mathbb{R}^d\ |\ a_{j\neq i}\geqslant 0\ \},\qquad i=1,\ldots,d.$$

\smallskip

\begin{defi}
A family $\cF$ of Metzler matrices is a {\rm product family} if there exist  compact sets $\cF_i\subset\cM_i,\ i=1,\ldots,d$, such that $\cF$ consists of all possible matrices with $i$-th row from $\cF_i$, for every $i=1,\ldots,d$. If for every $\cF_i$ we have $\cF_i\subset\mathbb{R}^d_+$, then $\cF$ is {\rm product family of non-negative matrices}. In either case, the sets $\cF_i$ are called {\rm uncertainty sets}.
\end{defi}

\nid The matrices belonging to the product families are constructed by independently taking $i$th row from the uncertainty set $\cF_i$. Moreover, product families can be topologically seen as product sets of their uncertainty sets: $\cF = \cF_1\times\cdots\times\cF_d$.\\

\bigskip

\begin{center}
\textbf{3.1. The greedy spectral method for the spectral abscissa}
\end{center}

\medskip

In \cite{Pro} \textit{spectral simplex method} for optimizing a spectral radius over non negative  product families was presented. There, it was also observed that the theory behind this procedure can be generalized for the optimization of the spectral abscissa over Metzler product families.\\ 

\cite{NePr} presents the \textit{greedy spectral simplex method} (or shorter \textit{greedy method}) for optimizing the spectral radius over positive product families, together with the modifications for the minimization on non-negative product families. Numerical experiments showed that this procedure is strikingly more efficient than the spectral simplex method. This efficiency was theoretically confirmed in \cite{PrCv}, where it was proved that greedy method has local quadratic convergence. Here we present a version of this algorithm for maximizing spectral abscissa over Metzler product families.\\

\newpage

\nid {\tt Algorithm 2: Maximizing spectral abscissa over Metzler product families}

\bigskip

\noindent {\tt Initialization.} Let $\cF = \cF_1\times\cdots\times\cF_d$ be a product family of Metzler matrices. Taking arbitrary $A^{(i)} \in \cF_i\, , \ i = 1, \ldots , d$, form a Metzler matrix $A_1 \in \cF$ with rows $A^{(1)}, \ldots , A^{(d)}$. Take its leading eigenvector $\bv_1$. Denote  $A_1^{(i)}\, = \, A^{(i)}, \ i = 1, \ldots , d$.\\

\smallskip

\noindent {\tt Main loop; the $k$th iteration.} We have a matrix $A_{k} \in \cF$ composed with rows $A_k^{(i)} \in \cF_i,
i = 1, \ldots , d$. Compute its leading eigenvector
$\bv_{k}$ (if it is not unique, take any of them) and for  $i = 1, \ldots d$, find a solution
$\hat{A}^{(i)} \in \cF_i$ of the problem
\begin{equation*}
\left\{
\begin{array}{l}
\langle A^{(i)} ,\bv_{k}\rangle \ \to \ \max\\
A^{(i)} \in \cF_i
\end{array}
\right.
\end{equation*} 
For each $i = 1, \ldots , d$, do:\\

\smallskip

If $\langle\hat{A}^{(i)}, \bv_{k}\rangle = \langle A_k^{(i)}, \bv_{k} \rangle$, then set $A_{k+1}^{(i)} = A_k^{(i)}$.

Otherwise, if  $\langle\hat{A}^{(i)}, \bv_{k}\rangle > \langle A_k^{(i)}, \bv_{k} \rangle$, set $A_{k+1}^{(i)} = \hat{A}^{(i)}$.\\
 
 \noindent Form the corresponding matrix~$A_{k+1}$.
 If the first case took place for all $i$, i.e. if 
 $A_{k+1} = A_k$,
 go to {\tt Termination}. Otherwise, go to $(k+1)$st iteration.\\
\smallskip 

\noindent {\tt Termination.} If $\bv_k > 0$, the procedure is finished; $A_k$ is maximal in each row\footnote{\ see Definition \ref{minimax}.} and $\eta_{\max} = \eta(A_k)$ is a solution. If $\bv_k$ has some zero components,
 then the family $\cF$ is reducible, and we need to stop the algorithm and 
 factorize $\cF$ (see section 3.4). {\hfill $\diamond$}\\

\smallskip

\nid \textbf{Remark.} The procedure for minimizing the spectral abscissa is exactly the same, except that we change the row if $\langle\hat{A}^{(i)}, \bv_{k}\rangle < \langle A_k^{(i)}, \bv_{k}\rangle$, and we omit the requirement that $\bv_k > 0$.

\bigskip

\begin{center}
\textbf{3.2. Theoretical results}
\end{center}

\medskip

For general Metzler matrices, the described greedy algorithm may  cycle. Hence, we one needs to implement some additional modifications. We will address this issue later on, but for now, we will focus on full Metzler families, since in this case the cycling does not occur. As with the positive and non-negative product families, the following theoretical results are the extension of the similar results from \cite{NePr2}.

\begin{lemma}\label{the_lemma}
Let $A$ be a Metzler matrix, $\bu\geqslant 0$ be a vector, and $\lambda\geqslant 0$ be a real number. Then $A\bu\geqslant\lambda\bu$ implies that $\eta(A)\geqslant\lambda$. If for a strictly positive vector $\bv$, we have $A\bv\leqslant\lambda\bv$, then $\eta(A)\leqslant\lambda$.
\end{lemma}
{\tt Proof.} Let $A\bu\geqslant\lambda\bu$. Since $A$ is Metzler, there exists $h\geqslant0$ such that $A+hI\geqslant 0$. We have $(A+hI)\bu\geqslant(\lambda+h)\bu$, and since the analogous of this Lemma for non-negative matrices and spectral radii holds \cite{NePr}, we obtain $\rho(A+hI)\geqslant\lambda+h$, and therefore $\eta(A)\geqslant\lambda$ (Lemma \ref{trans}). The proof of the second statement follows by the same reasoning.
{\hfill $\Box$}

\begin{cor}\label{the_cor} Let $A$ be a Metzler matrix, $\bu\geqslant 0$ a vector, and $\lambda\geqslant 0$ a real number. Then $A\bu>\lambda\bu$ implies that $\eta(A)>\lambda$. If for a strictly positive vector $\bv$, we have $A\bv<\lambda\bv$, then $\eta(A)<\lambda$.
\end{cor}
{\tt Proof.} Let $A\bu>\lambda\bu$. The statement $\eta(A)<\lambda$ is in direct contradiction with Lemma \ref{the_lemma} and $\eta(A)=\lambda$ would imply that $\lambda$ is an eigenvalue of matrix $A$, which is not true. Hence, $\eta(A)>\lambda$ has to hold. The proof of the second statement follows the same reasoning.
{\hfill $\Box$}\\

\begin{defi}\label{minimax}
We say that the matrix $A$ is {\rm maximal in each row with the respect to} $\bv$ if $\langle A^{(i)},\bv\rangle = \max\limits_{\bx\in\cF_i}\langle\bx,\bv\rangle$ for all $i=1,\ldots,d$. Similarly we define minimality in each row.
\end{defi}

\nid The following Proposition stems directly from Lemma \ref{the_lemma}.

\begin{prop}\label{the_prop}
Let matrix $A$ belong to a Metzler product family $\cF$ and $\bv\geqslant 0$ be its leading eigenvector. Then\\
1) if $\bv>0$ and $A$ is maximal in each row with the respect to $\bv$, then $\eta(A) = \max\limits_{X\in\cF}\eta(X)$.\\
2) if $A$ is minimal in each row with the respect to $\bv$, then $\eta(A) = \min\limits_{X\in\cF}\eta(X)$.
{\hfill $\Box$}
\end{prop}

\begin{prop}\label{the_propII}
If a Metzler matrix $A\in\cF$ is full, then it has the maximal spectral abscissa in $\cF$ if and only if $A$ is maximal in each row with the respect to its (unique) leading eigenvector. The same is true for minimization.
\end{prop}
{\tt Proof.} {\tt Maximization.} Let $A$ be maximal in each row w.t.r. to $\bv$. Since for a full Metzler matrix its leading eigenvector is strictly positive, we can apply $1)$ of the previous Proposition, and therefore the maximal spectral abscissa is indeed achieved for matrix $A$.\\

Conversely, assume that $A$ maximizes the spectral abscissa on $\cF$, but it is not maximal w.t.r. to its leading eigenvector $\bv$ in each row. Construct matrix $A'\in\cF$ by changing all the necessary rows of the matrix $A$, so that $A'$ becomes maximal w.t.r. to $\bv$ in each row. Thus we have $A'\bv>A\bv=\eta(A)\bv$. This, in turn, results in $\eta(A')>\eta(A)$ (Corollary \ref{the_cor}), which is incorrect. So, $A$ has to be maximal in each row w.t.r. to its eigenvector.\\

\noindent {\tt Minimization.} Suppose now that $A$ is minimal in each row w.t.r. to $\bv$. One direction of the equivalence is basically $2)$ of Proposition \ref{the_prop}. The proof for the other direction goes the similar way as with the maximization, taking into account the  strict positivity of the vector $\bv$. 
{\hfill $\Box$}

\begin{prop}\label{the_propIII} 
For every Metzler product family there exists a matrix which is maximal (minimal) in each row with the respect to the one of its leading eigenvalues.
\end{prop}

\nid {\tt Proof.} Let $\varepsilon'>0$ be such that the shifted product family $\cF_{\varepsilon'} = \cF + \varepsilon' E$ is full Metzler family, and let $A_{\varepsilon'}\in\cF_{\varepsilon'}$ be the matrix with maximal spectral radius. Since this matrix is full, by Proposition \ref{the_propII} it has to be maximal in each row. For all $0<\varepsilon\leqslant\varepsilon'$, we associate one such a matrix $A_\varepsilon$. By compactness, there exists a sequence $\{\varepsilon_k\}_{k\in\mathbb{N}}$ such that $\varepsilon_k\rightarrow 0$ as $k\rightarrow +\infty$. Hence, matrices $A_{\varepsilon_k}$ converges to a matrix $A\in\cF$ and their respective leading eigenvectors converge to a nonzero vector $\bv$. By continuity, $\bv$ is a leading eigenvector of $A$, and $A$ is maximal in each row w.t.r. to $\bv$.
{\hfill $\Box$}

\begin{theorem}\label{the_them}
If $\cF$ is a product family of full Metzler matrices, the greedy spectral simplex method terminates in finite time. 
\end{theorem}

\noindent {\tt Proof.} Assume we use the greedy method for maximizing spectral abscissa (the same goes for minimization). Proposition \ref{the_propIII} guarantees the existence of a matrix maximal in each row, and Proposition \ref{the_propII} that this matrix is the optimal one. Starting from matrix $A_1\in\cF$ and iterating trough the algorithm we obtain the sequence of matrices $A_2,A_3,\ldots$ with the sequence of their respective spectral abscissas $\eta(A_2),\eta(A_3),\ldots$, which is increasing. Moreover, each $i$th row of of the given matrices is some vertex of the polyhedron $\cF_i$, so the number of total states is finite. Therefore, we arrive at our solution in finite number of iterations, increasing the spectral abscissa in each one, until we reach the optimal matrix.
{\hfill $\Box$}

\newpage

\begin{center}
\textbf{3.3. Selective greedy method for spectral abscissa}
\end{center}

\medskip

To resolve the issue of cycling in \cite{PrCv} the \textit{selective greedy method} was proposed. This method is as efficient as the greedy method, and it does not cycle even if the matrices are very sparse.\\

\smallskip

\begin{defi}
The {\rm selected leading eigenvector} of a Metzler matrix $A$ is the limit $\lim\limits_{\varepsilon \to 0}\bv_{\varepsilon}$, 
where $\bv_{\varepsilon}$ is the normalized leading eigenvector of the 
perturbed matrix~$A_{\varepsilon} \, = \, A\, +\, \varepsilon E$. 
\end{defi}

\noindent Notice that in this definition we have extended the notion of the selected leading eigenvector from \cite{PrCv} to include Metzler matrices as well, which will be fully justified further in the text.\\

\smallskip

\begin{prop}\label{pow.mthd}
The power method $\bx_{k+1} = A\bx_k, \, k \ge 0,$ where $A$ is a non-negative matrix, applied to the  initial vector $\bx_0 = \be$ converges to the selected leading eigenvector.
{\hfill $\Box$}
\end{prop}

\begin{defi}\label{selecta}
The greedy method with the selected leading eigenvectors~$\bv_k$ 
in all iterations is called {\rm selective greedy method}.
\end{defi}

\begin{theorem}\label{selecta.th}
The selective greedy method, applied to non-negative matrices, does not cycle.
{\hfill $\Box$}
\end{theorem}

\smallskip

Using the selective greedy method on Metzler matrices to avoid cycling in the sparse case seems like a logical move. However, it may happen that the power method applied to a strict Metzler matrix, starting with the initial vector $\bv_0 = \be$, does not converge.\\ 

\smallskip

\noindent \textbf{Example.} We apply the power method to the matrix 
$$
A = 
\left(
\begin{array}{rrr}
-2  & \phantom{-}2 & \phantom{-}0\\
\phantom{-}0 & -6 & \phantom{-}5\\
\phantom{-}2 & \phantom{-}2 & -9\\
\end{array}
\right).
$$
\nid Starting from the vector $\bv_0=\be$, and obtain the following sequence of vectors:\\

\nid $\bv_1 = (0.242,\ 0,\ -0.97)\\
\bv_2 = (-0.025,\ -0.507,\ 0.862)\\
\bv_3 = (-0.094,\ -0.649,\ -0.755)\\
\bv_4 = (0.138,\ -0.694,\ 0.707)\\
\cdots$\\

\noindent It is clear that the sequence of vectors $\bv_k$ is not going to converge, because the sign pattern keeps changing with every iteration.\\

\smallskip

Having in mind that a Metzler matrix can be translated to a corresponding non-negative matrix, we propose a \textit{translative power method}.

\begin{defi}
Let $A$ be a Metzler matrix. A power method applied to a matrix $\tilde{A}=A+hI$, where $h\geqslant0$ is the minimal number for which $\tilde{A}$ is non-negative is called {\rm translative power method}.
\end{defi}

\begin{prop}\label{trans.power}
The translative power method $\bx_{k+1} = A\bx_k, \, k \ge 0,$ where $A$ is Metzler matrix, applied to the initial vector $\bx_0 = \be$ converges to the selected leading eigenvector. 
\end{prop}

\noindent {\tt Proof.} If $A$ is non-negative, then translative power method is just the regular power method and we can use Proposition \ref{pow.mthd}. Suppose now that $A$ is strictly Metzler. Applying a translative power method on it, taking $\be$ as an initial vector, we obtain a selected leading eigenvector $\bv$ for a non-negative matrix $\tilde{A} = A + hI$. This eigenvector $\bv = \lim\limits_{\varepsilon \to 0}\bv_{\varepsilon}$ is a limit of\vspace*{0.1cm} a sequence of leading eigenvectors of perturbed positive matrices $\tilde{A}_\varepsilon = \tilde{A} +\varepsilon E$. Since the leading eigenvectors $\bv_\varepsilon$ for matrices $\tilde{A}_\varepsilon$ are also the leading eigenvectors for the perturbed Metzler matrices $A_\varepsilon = A + \varepsilon E$, we have, by the definition, that $\bv$ is a selected leading eigenvector of the matrix $A$.
{\hfill $\Box$}\\

\smallskip

\noindent Proposition \ref{trans.power} allows us to use the Definition \ref{selecta} when dealing with Metzler matrices as well, having in mind that in this case we need to resort to translative power method. Moreover,

\begin{theorem}
The selective greedy method, applied to Metzler matrices, does not cycle.
\end{theorem}
\noindent {\tt Proof.} Assume we have a product family $\cF$ that contains strict Metzler matrices, since the case of non-negative families is already resolved by Theorem \ref{selecta.th}. Cycling of a selective greedy method on this family would also imply the cycling of the greedy spectral simplex method on the family $\cF + \varepsilon E$ of full Metzler matrices, which is, given Theorem \ref{the_them}, impossible.
{\hfill $\Box$}\\

\noindent \textbf{Remark.} As discussed in \cite{PrCv}, when implementing selective greedy method we compute the selected leading eigenvectors by applying the power method not on the non-negative matrices $X_k$, obtained iterating through the procedure, but on the matrices $X_k + I$. This is due to the fact that, for any $k$, both matrices $X_k$ and $X_k + I$ have the same selected leading eigenvector, while for the matrix $X_k + I$ it will be unique, and maybe multiple. The uniqueness of this eigenvector is a guarantee that the power method will converge.\\
For the same reasons, if we obtain a strict Metzler matrix $X$ in some iteration, we will actually compute the selected leading eigenvector by applying power method on a non-negative matrix $X + (h+1)I$, where $h = |\min_i x_{ii}|$.\\
In addition, all remarks regarding the cycling due to the computational errors and how to avoid it, given in \cite{PrCv}, apply for the Metzler matrices as well.\\

\newpage

\begin{center}
\textbf{3.4. Optimizing spectral abscissa for reducible families}
\end{center}

\medskip

If we finish the (selective) greedy procedure with a leading eigenvector $\bv_k$ that has zero entries, we might not obtain an optimal solution. In this case the product family $\cF$ is reducible. This means that the linear space spanned by the vectors with the same support as $\bv_k$ is invariant with the respect to all matrices of $\cF$. One way to resolve this is to use the Frobenius factorization and run the maximization procedure on irreducible blocks $\cF^{(j)}$, constructing the optimal solution as the block matrix. The optimal value $\eta_\max$ is the largest of all optimal values $\eta_\max^{(j)}$ among all blocks. The algorithm for the Frobenius factorization can be found in \cite{Tar}.\\

In practice, however, the Frobenius factorisation usually takes more time than the run of the greedy algorithm for the whole family. That is why it makes sense to avoid the factorisation procedure by making the family irreducible. 
Since having only one irreducible matrix is enough to make the whole family irreducible, a simple strategy  is to simply include a matrix $H = \alpha P - \beta I$ to the family $\cF$, where $P$ is a cyclic permutation matrix, $\alpha > 0$, and $\beta\geqslant 0$. The matrix $P$ is irreducible, and so the family $\cF\cup\{H\}$ is irreducible as well, hence the greedy method has to finish with a positive leading eigenvector. If the final matrix has no rows from $H$, then it is optimal for 
the original family $\cF$. Otherwise we have to choose another pair $\alpha, \beta$, setting them to be smaller. However, if trying out various pairs of $\alpha, \beta$ produces the final matrix with some rows from $H$, we need to resort to Frobenius factorization.\\  

\bigskip

\begin{center}
\textbf{3.5. Numerical results}
\end{center}

\medskip

We test the selective greedy method on both full Metzler product families and sparse Metzler product families with density parameters $\gamma_i$ (the percentage of nonzero entries of the vectors belonging to the uncertainty set $\mathcal{F}_i$).\\

\nid The results of tests on full Metzler product families, as dimension $d$ and the size of the uncertainty sets $N$ vary, are given in Tables 1 and 2, for maximization and minimization, respectively. In Tables 3 and 4 we report the behaviour of selective greedy algorithm on sparse Metzler families. For each uncertainty set $\mathcal{F}_i$ the density parameter $\gamma_i$ is randomly chosen from the interval $9-15\%$ and the elements of the set $\mathcal{F}_i$ are randomly generated in accordance to it.\\

\begin{table}[H]
\begin{center}
\begin{tabular}{c|c c c c}
$N\setminus d$ & 25  & 100 & 500 & 2000\\
\hline
50 & 3 & 3 & 3.1 & 3\\
100 & 3.2 & 3.1 & 3.1 & 3.1\\
250 & 3.2 & 3.2 & 3.1 & 3.1\\
\end{tabular}
\caption*{{\footnotesize Table 1a: Average number of iterations for maximization, for full Metzler families}}
\end{center}
\end{table}

\begin{table}[H]
\begin{center}
\begin{tabular}{c|c c c c}
$N\setminus d$ & 25  & 100 & 500 & 2000\\
\hline
50 & 0.02s & 0.09s & 0.5s & 2.77s\\
100 & 0.04s & 0.19s & 1.02s & 6.5s\\
250 & 0.1s & 0.44s & 2.19s & 136.77s\\
\end{tabular}
\caption*{{\footnotesize Table 1b: Average computing time for maximization, for full Metzler families}}
\end{center}
\end{table}

\begin{table}[H]
\begin{center}
\begin{tabular}{c|c c c c}
$N\setminus d$ & 25  & 100 & 500 & 2000\\
\hline
50 & 3.1 & 3.1 & 3.1 & 3.1\\
100 & 3.1 & 3 & 3 & 3.1\\
250 & 3.3 & 3.2 & 3.1 & 3.1\\
\end{tabular}
\caption*{{\footnotesize Table 2a: Average number of iterations for minimization, for full Metzler families}}
\end{center}
\end{table}

\begin{table}[H]
\begin{center}
\begin{tabular}{c|c c c c}
$N\setminus d$ & 25  & 100 & 500 & 2000\\
\hline
50 & 0.02s & 0.09s & 0.54s & 2.5s\\
100 & 0.04s & 0.19s & 1s & 4.6s\\
250 & 0.1s & 0.45s & 2.32s & 129.35s\\
\end{tabular}
\caption*{{\footnotesize Table 2b: Average computing time for minimization, for full Metzler families}}
\end{center}
\end{table}
  
\begin{table}[H]
\begin{center}
\begin{tabular}{c|c c c c}
$N\setminus d$ & 25  & 100 & 500 & 2000\\
\hline
50 & 5.4 & 4.4 & 4.2 & 4\\
100 & 5.7 & 4.8 & 4.2 & 4.2 \\
250 & 6.1 & 4.8 & 4.2 & 4.2\\
\end{tabular}
\caption*{{\footnotesize Table 3a: Average number of iterations for maximization, for sparse Metzler families}}
\end{center}
\end{table}

\begin{table}[H]
\begin{center}
\begin{tabular}{c|c c c c}
$N\setminus d$ & 25  & 100 & 500 & 2000\\
\hline
50 & 0.04s & 0.13s & 0.64s & 3.16s\\
100 & 0.07s & 0.26s & 1.3s & 6.03s\\
250 & 0.19s & 0.69s & 3.13s & 151.39s\\
\end{tabular}
\caption*{{\footnotesize Table 3b: Average computing time for maximization, for sparse Metzler families}}
\end{center}
\end{table}

\begin{table}[H]
\begin{center}
\begin{tabular}{c|c c c c}
$N\setminus d$ & 25  & 100 & 500 & 2000\\
\hline
50 & 7.2 & 5 & 4.1 & 4.1\\
100 & 6.9 & 5.2 & 4.3 & 4.1 \\
250 & 12.2 & 5.4 & 4.6 & 4.2 \\
\end{tabular}
\caption*{{\footnotesize Table 4a: Average number of iterations for minimization, for sparse Metzler families}}
\end{center}
\end{table}

\begin{table}[H]
\begin{center}
\begin{tabular}{c|c c c c}
$N\setminus d$ & 25  & 100 & 500 & 2000\\
\hline
50 & 0.87s & 0.14s & 0.66s & 3.15s\\
100 & 0.14s & 0.28s & 1.29s & 6.9s\\
250 & 0.46s & 0.85s & 3.35s & 169.89s\\
\end{tabular}
\caption*{{\footnotesize Table 4b: Average computing time for minimization, for sparse Metzler families}}
\end{center}
\end{table}

\nid Table 5 shows how the number of iterations and computing time vary as the density parameter is changed. The dimension is kept fixed at $d = 600$ and the cardinality of each product set at $\vert\mathcal{F}_i\vert=200$, while we vary the interval from which $\gamma_i$ takes value.\\

\begin{table}[H]
\begin{center}
\begin{tabular}{c|c c c c}
$\gamma_i$ & 9-15  & 16-21 & 22-51 & 52-76\\
\hline
MAX & 4 & 4 & 4 & 3.9\\
MIN & 4.2 & 4 & 4.2 & 3.6\\
\end{tabular}
\caption*{{\footnotesize Table 5a: Effects of sparsity, number of iterations}}
\end{center}
\end{table}

\begin{table}[H]
\begin{center}
\begin{tabular}{c|c c c c}
$\gamma_i$ & 9-15  & 16-21 & 22-51 & 52-76\\
\hline
MAX & 2.54s & 2.5s & 2.5s & 2.49s\\
MIN & 2.76s & 2.66s & 2.96s & 2.33s\\
\end{tabular}
\caption*{{\footnotesize Table 5b: Effects of sparsity, computing time}}
\end{center}
\end{table}

\bigskip

\begin{center}\label{sec4}
\large{\textbf{4. Problems of (de)stabilization in $l_\infty$ norm}}
\end{center}

\smallskip

Let $A$ be a Metzler matrix, and let
$$\cB_\varepsilon(A) = \{X\in\cM\ |\ \|A-X\|\leqslant\varepsilon\}$$
\nid denote an $\varepsilon$-ball of Metzler matrices around $A$. Optimization of the spectral abscissa on $\cB_\varepsilon(A)$ is a crucial tool for the Hurwitz (de)stabilization. The greedy method developed in Section 3 can be applied on $\cB_\varepsilon$, since it can be considered as a product set of the balls
$$\cB_\varepsilon^{(i)} = \{X^{(i)}\in\cM_i\ |\ \|A^{(i)}-X^{(i)}\|\leqslant\varepsilon\},\qquad i=1,\ldots,d.$$

\medskip

\begin{center}
\textbf{4.1. Closest Hurwitz unstable Metzler matrix in the $l_\infty$ norm}
\end{center}

\medskip

Let us first consider the problem of Schur destabilization of non-negative matrix in $l_\infty$ norm: if $A$ is a non-negative matrix with $\rho(A) < 1$, find a matrix $X$ that satisfies
\begin{equation}\label{schur}
\left\{
\begin{array}{l}
\|X - A\|\ \rightarrow\min\\
\rho(X) = 1.
\end{array}
\right.
\end{equation}
An explicit solution to (\ref{schur}) was given in \cite{NePr}. An important point is that the solution is always a non-negative matrix, entrywise bigger than $A$. 

\noindent\begin{theorem}{\rm \cite{NePr}}\label{schunstable}
The optimal value $\tau_*$ of the problem {\rm (\ref{schur})} is the reciprocal of the biggest component of the vector $(I - A)^{-1}\be$. Let $k$ be the index of that component. Then the optimal solution is the matrix
$$X = A + \tau_*E_k.$$
{\hfill $\Box$}
\end{theorem}

\smallskip

We can also solve a more general problem: for a given $h > 0$ and non-negative matrix $A$ with $\rho(A) < h$, find the closest matrix $X$ having $\rho(X) = h$, i.e., find a solution to
\begin{equation}\label{hschurmain}
\left\{
\begin{array}{l}
\|X - A\|\ \rightarrow\min\\
\rho(X) = h.
\end{array}
\right.
\end{equation}
Here, as with the $h = 1$ case, the solution $X$ is also non-negative and entrywise bigger than $A$. We can also prove a generalization of Theorem \ref{schunstable}:

\noindent\begin{theorem}\label{hschunstable}
The optimal value $\tau_*$ of the problem {\rm (\ref{hschurmain})} is the reciprocal of the biggest component of the vector $(hI - A)^{-1}\be$. Let $k$ be the index of that component. Then the optimal solution is the matrix
\begin{equation}\label{htauschur}
X = A + \tau_*E_k.
\end{equation}
\end{theorem}

\smallskip

\noindent {\tt Proof.} The optimal matrix $X_*$ for (\ref{hschurmain}) is also a solution to the maximization problem
\begin{equation*}
\left\{
\begin{array}{l}
\rho(X)\ \rightarrow\max\\
\|X - A\|\leqslant\tau,\ X\geqslant A.
\end{array}
\right.
\end{equation*}
having the optimal value $\tau = \tau_*$. Let us characterize this matrix for arbitrary $\tau$. From  Proposition \ref{the_prop} and Proposition \ref{the_propIII}, applied to the non-negative product families and spectral radii, we can conclude that $X$ is maximal in each row for the product family with the uncertainty sets:
\begin{equation*}
\left.
\begin{array}{lll}
\cB_\tau^+(A) & = & \cB_\tau(A)\cap\{X\geqslant 0\ |X\geqslant A\}\\
\\
\ & = & \{X\geqslant 0\ |X\geqslant A, \langle(X-A)\be,\be_i\rangle\leqslant\tau,\ i=1,\ldots,d\}.
\end{array}
\right.
\end{equation*}
\nid Conversely, every matrix $X$ with $\rho(X) = h$ which is maximal in each row w.r.t. a strictly positive leading eigenvector, solves  (\ref{hschurmain}).\\

\nid Any matrix $X\in\cB_\tau^+(A)$ with the leading eigenvector $\bv$ is optimal in the $i$th row if and only if the scalar product $\langle X^{(i)} - A^{(i)},\bv\rangle$ is maximal under the constraint $\langle X^{(i)} - A^{(i)},\be\rangle = \tau$. This maximum is equal to $r\tau$, where $r$ is the maximal component of the leading eigenvector $\bv$. Denote the index of this component by $k$. Then
$$X^{(i)} - A^{(i)} = \tau\be_k,\qquad i=1,\ldots,d.$$
\nid Hence, if $X$ is maximal in each row, we have 
$$X = A + \tau\be\be^T_k = A + \tau E_k.$$
\nid Furthermore, since each set $\cB_\tau^+(A^{(i)})$ contains a strictly positive point, we have 
$$v_i = (X\bv)_i = \max\limits_{\bx\in\cB_\tau(A^{(i)})}\langle\bx,\bv\rangle > 0.$$
\nid Therefore, the leading eigenvector $\bv$ is strictly positive so, by Proposition \ref{the_prop}, the matrix $X$ maximizes the spectral radius on the product family $\cB_\tau^+(A)$.\\

Thus, the optimal matrix has the form (\ref{htauschur}) for some $k$, and $\|X-A\|=\tau_*$. It remains to find $k$ for which the value of $\tau$ is minimal.\\

\nid Since $\rho(A + \tau_*E_k) = h$, it follows that $\tau_*$ is the smallest positive root of the equation
\begin{equation}\label{eq.1}
\det(A - hI + \tau E_k) = 0.
\end{equation}
Since $\rho(\frac{A}{h}) < 1$ we have $(hI - A)^{-1} = \frac{1}{h}(I - \frac{A}{h})^{-1} = \frac{1}{h}\sum_{j=0}^{\infty}(\frac{A}{h})^j\geqslant 0$. Multiplying equation (\ref{eq.1}) by $\det(-(hI - A)^{-1})$ we obtain
$$\det(I - \tau(hI - A)^{-1}E_k) = 0.$$
The matrix $\tau(hI - A)^{-1}E_k$ has only one nonzero column. This is the $k$th column, equal to $\tau(hI - A)^{-1}\be$. Hence,
$$\det(I - \tau(hI - A)^{-1}E_k) = 1 - \tau[(hI - A)^{-1}\be]_k.$$
We conclude that the minimal $\tau$ corresponds to the biggest component of this vector. {\hfill $\Box$}\\

We can now move to a Hurwitz destabilization, or more formally put: for a given Metzler matrix $A$ with $\eta(A) < 0$, find a solution to
\begin{equation}\label{hurwitz}
\left\{
\begin{array}{l}
\|X - A\|\ \rightarrow\min\\
\eta(X) = 0.
\end{array}
\right.
\end{equation}

\noindent\begin{lemma}\label{lemma.hur}
Let a matrix $X$ be a solution to the problem (\ref{hurwitz}). Then, $X$ is Metzler and $X \geqslant A$.
\end{lemma}

\noindent {\tt Proof.} First, let us show that $X$ is Metzler. Assume the contrary, that it has some negative off-diagonal entries. We consider a matrix $X'$ with entries
\begin{equation*}
x'_{ij} = 
\left\{
\begin{array}{cc}
x_{ij}, & x_{ij}\geqslant 0\ {\rm or}\ i = j \\
|x_{ij}|, & x_{ij} < 0\ {\rm and}\ i\neq j.
\end{array}
\right.
\end{equation*}
\nid We have $\|X' - A\| < \|X - A\|$ and $\eta(X') < 0$. Since $X'$ is Metzler, there exists $h > 0$ such that $X'_h = X' + hI$ is non-negative. In addition, from Lemma \ref{trans}, we obtain $\rho(X'_h) < h$.\\ 

\nid Now, define the matrix $X_h = X + hI$ which contains off-diagonal negative entries. From the inequality $\|X_h^k\|\ < \|(X'_h)^k\|$, using Gelfand's formula for spectral radius \cite{HoJo}, we arrive at $\rho(X_h) < \rho(X'_h)$. Since the largest real part of all the eigenvalues of $X$ is equal to zero, we have
$$ 0 = \eta(X)\leqslant\rho(X). $$

\nid Adding $h$ on both sides of inequality and using $\rho(X) + h = \rho(X + hI)$, we obtain $h < \rho(X'_h)$, which is a contradiction. Therefore, matrix $X$ must not contain negative off-diagonal entries, i.e. it has to be Metzler.\\  

Now let us prove that $X \geqslant A$. Assume the converse, that there exist some entries $x_{ij}$ of $X$ such that $x_{ij} < a_{ij}$. Define the matrix $\tilde{X}$ with entries
\begin{equation*}
\tilde{x}_{ij} = 
\left\{
\begin{array}{cc}
x_{ij}, & x_{ij}\geqslant a_{ij}\\
a_{ij}, & x_{ij} < a_{ij}
\end{array}
\right.
\end{equation*}
\nid $\tilde{X}$ is Metzler, $\tilde{X}\geqslant A$ and $\tilde{X}\geqslant X$. By the monotonicity of spectral abscissa we have $\eta(\tilde{X})\geqslant\eta(X)$. $X$ is the closest Hurwitz unstable Metzler matrix, so $\|\tilde{X} - A\|\ \geqslant\|X - A\|$. This is impossible, since $\|\tilde{X} - A\|\ < \|X - A\|$ holds. Hence, $X$ has to be entrywise bigger than $A$. {\hfill $\Box$}

\noindent\begin{theorem}
The optimal value $\tau_*$ of the problem {\rm (\ref{hurwitz})} is the reciprocal of the biggest component of the vector $-A^{-1}\be$. Let $k$ be the index of that component. Then the optimal solution is the matrix
\begin{equation}\label{eq.2}
X = A + \tau_*E_k.
\end{equation}
\end{theorem}

\medskip

\noindent {\tt Proof.} Lemma \ref{lem2} ensures invertibility of the matrix $A$. We also remark that a Hurwitz stable Metzler matrix is necessary a strict Metzler matrix: if it were non-negative, than the biggest real part of its eigenvalues would also be non-negative, which is not the case.\\

\nid Now, let $h > 0$ be such that the matrix $\tilde{A} = A + hI$ is non-negative. We have $\rho(\tilde{A}) < h$. Using Theorem \ref{hschunstable} we obtain the closest non-negative matrix to $\tilde{A}$. Denote it by  $\tilde{X}$. We have $\rho(\tilde{X}) = h$ and 
$$\tilde{X} = \tilde{A} + \tau_*E_k,$$
\nid where $\tau_*$ is the reciprocal of the biggest component of the vector $(hI - \tilde{A})^{-1}\be = -A^{-1}\be$. Let $X = \tilde{X} - hI$. $X$ is Metzler with $\eta(X) = 0$, and satisfies (\ref{eq.2}).\\

Now we need to check that $X$ is really the closest Hurwitz unstable Metzler matrix for $A$. Assume that $Y \neq X$ is the true solution of (\ref{hurwitz}). We have $\|Y-A\|\ < \|X-A\| = \tau_*$. Define $\tilde{Y} = Y + hI$. $\tilde{Y}$ is entrywise bigger than $\tilde{A}$, and $\rho(\tilde{Y}) = h$. Since $\tilde{X}$ is the closest non-negative matrix to $\tilde{A}$ having spectral radius equal to $h$, the following is true:
$$\|Y-A\|\ = \|\tilde{Y} - \tilde{A}\|\ \geqslant \|\tilde{X} - \tilde{A}\|\ = \tau_*,$$
which is a contradiction. Therefore, our $X$ is really the optimal solution.    
{\hfill $\Box$}\\

\medskip

\noindent{\bf Example.}
For a given Hurwitz stable Metzler matrix

$$
A = 
\left(
\begin{array}{rrrrr}
-4  & {\ } 0 & {\ } 0 & {\ } 0 & {\ } 4 {\ }\\
{\ } 0  & -2 & {\ } 0 & {\ } 2  & {\ } 0 {\ }\\
{\ } 0  & {\ } 2 & -1 & {\ } 0 & {\ } 0 {\ }\\
{\ } 0  & {\ } 0 & {\ } 0 & -4 & {\ } 0 {\ }\\
{\ } 0  & {\ } 0 & {\ } 0 & {\ } 3 & -9 {\ }\\
\end{array} 
\right)
$$\\
with $\eta(A) = -1$, the biggest component of the vector $-A^{-1}\be$ is the third one, and it is equal to $2.5$. Thus we get the matrix
$$
X = 
\left(
\begin{array}{rrcrr}
-4  & {\ } 0 & {\ }{\ } 0.4 & {\ } 0 & {\ } 4 {\ }\\
{\ } 0  & -2 & {\ }{\ } 0.4 & {\ } 2  & {\ } 0 {\ }\\
{\ } 0  & {\ } 2 & -0.6 & {\ } 0 & {\ } 0 {\ }\\
{\ } 0  & {\ } 0 & {\ }{\ } 0.4 & -4 & {\ } 0 {\ }\\
{\ } 0  & {\ } 0 & {\ }{\ } 0.4 & {\ } 3 & -9 {\ }\\
\end{array} 
\right)
$$\\
as the closest Hurwitz unstable.\\

\bigskip

\begin{center}
\textbf{4.2. Closest Hurwitz stable Metzler matrix in the $l_\infty$ norm}
\end{center}

\medskip

The problem of Schur stabilization for non-negative matrices in $l_\infty$ norm was also considered in \cite{NePr}: if $A$ is a given non-negative matrix with $\rho(A) > 1$, find a non-negative matrix $X$ that satisfies
\begin{equation}\label{stab.schur}
\left\{
\begin{array}{l}
\|X - A\|\ \rightarrow\min:\ X\geqslant 0,\\
\rho(X) = 1.
\end{array}
\right.
\end{equation}
\nid Notice that, in contrast to the problem of Schur destabilization (\ref{schur}), here we have to impose the non-negativity condition on our solution.\\

The story is the same with the Hurwitz (de)stabilization of Metzler matrix. For Hurwitz destabilization (\ref{hurwitz}), requesting for the solution to be Metzler is redundant, a fact confirmed by Lemma \ref{lemma.hur}. However, this is not the case with the Hurwitz stabilization; here, we have to explicitly impose the restriction on our solution to be Metzler. So, our problem can be written in the following manner:
\begin{equation}\label{stab.hurwitz}
\left\{
\begin{array}{l}
\|X - A\|\ \rightarrow\min:\ X\in\cM,\\
\eta(X) = 0.
\end{array}
\right.
\end{equation}
\nid To obtain the solution, we solve the related problem
\begin{equation*}
\left\{
\begin{array}{l}
\eta(X)\rightarrow\min:\ \|X - A\|\ \leqslant\tau,\\
A\geqslant X,\ X\in\cM.
\end{array}
\right.
\end{equation*}
\nid and use the strategy that implements the bisection in $\tau$, together with the greedy procedure for minimizing the spectral abscissa described in Section 3.\\

 First, we describe how to compute matrix $X$ explicitly in each iteration of the Algorithm 2 for minimization, thus finding the optimal solution on the ball of radius $\tau$.\\

\nid Fix some $\tau > 0$, and let $X'$ be a matrix obtained in some iteration with $\bv\geqslant 0$ as its leading eigenvector. We rearrange positive entries of $\bv$: $v_{j_1}\geqslant\cdots\geqslant v_{j_m}$, where $\cS$ is support of $\bv$, and $|\cS|\ = m$. In the next iteration of the greedy procedure we construct the matrix $X = (x_{ij})$:\\
\\
For each $i\in\cS$ we solve
$$\langle\X^{(i)},\bv\rangle\rightarrow\min,\qquad X^{(i)}\in\cF_i$$
which we can rewrite as
\begin{equation}\label{eq.3}
\sum_{k=1}^m x_{ij_k}v_{j_k}\rightarrow\min :
\left\{
\begin{array}{l}
\sum_{k=1}^m x_{ij_k}\geqslant -\tau + \sum_{k=1}^m a_{ij_k},\\
x_{ij_k}\geqslant 0,\ i\neq j_k.
\end{array}
\right.
\end{equation}
\nid We can solve (\ref{eq.3}) explicitly:
\begin{equation}\label{expl}
x_{ij_k} =
\left\{
\begin{array}{lr}
0, & k < l_i\\
-\tau + \sum_{s=1}^{l_i} a_{ij_s}, & k = l_i\\
a_{ij_k}, & k > l_i
\end{array}
\right.
\end{equation}
\nid where
\begin{equation}\label{elle}
l_i = \min \Big\{i,\ \min \{l\in\cS|\sum_{s=1}^{l} a_{ij_s} > \tau \}\Big\}.
\end{equation}
\nid For $i\in\cS$, but $j\notin\cS$, we take $x_{ij} = a_{ij}$; and if $i\notin\cS$ we put $X^{(i)} = X^{'(i)}$.\\

\medskip

\nid {\tt Alg. 3: Computing the closest Hurwitz stable Metzler matrix in $l_\infty$ norm}

\bigskip

\nid {\tt Step 0.} Take $\frac{\|A\|}{2}$ as the starting value for $\tau$.\\

\noindent {\tt Step 1.} Start the selective greedy method for minimizing the spectral abscissa on the ball of Metzler matrices $\cB_\tau(A)$ (applying (\ref{expl})). Iterate until a matrix $X$ is obtained with $\eta(X) < 0$. When this is done, stop the greedy procedure, compute its leading eigenvector $\bv$, and proceed to the next step.\\

\nid If greedy procedure finishes finding the matrix with non-negative minimal spectral abscissa on $\cB_\tau(A)$, keep implementing the bisection on $\tau$. Do this until the ball $\cB_\tau(A)$ which contains a matrix with negative spectral abscissa is obtained.\\

\noindent {\tt Step 2.} Construct matrices $C = (c_{ij})$ and $R$, as follows. For each $i\in\cS$
\begin{equation*}
c_{ij_k} =
\left\{
\begin{array}{lr}
0, & k < l_i\\
\sum_{s=1}^{l_i} a_{ij_s}, & k = l_i\\
a_{ij_k}, & k > l_i
\end{array}
\right.
\end{equation*}
\nid where $l_i$ is given by (\ref{elle}), while $R$ is a boolean matrix having ones on positions $(i,l_i)$ and zeros in all other places. If some of the indices $i,j$ are not in the support, then we put $c_{ij} = x_{ij}$. We can write $X = C - \tau R$.\\

\nid Denote by $\tau_*$ a potential optimal value of the problem (\ref{stab.hurwitz}). To determine this value we proceed to the next step.\\

\noindent {\tt Step 3.} Since $\eta(X) < 0$, we have by Lemma \ref{lem2} that $X$ is invertible and
$$-X^{-1} = - (C - \tau R)^{-1} \geqslant 0.$$
Since $0 = \eta(C - \tau_* R) = \eta(C - \tau R + (\tau - \tau_*)R)$, it follows that $\det(-(C - \tau R) - (\tau - \tau_*)R) = 0$. From here we have\\
\begin{equation*}
\det\Big(\frac{1}{\tau - \tau_*}I + (C-\tau R)^{-1}R\Big) = 0.
\end{equation*}\\
\nid Matrix $-(C - \tau R)^{-1}R$  is non-negative and $\frac{1}{\tau - \tau_*}$ is its (positive) leading eigenvalue. Therefore, in this step we find the potential optimal value $\tau_*$ by computing the leading eigenvalue $\lambda$ of the matrix $-(C - \tau R)^{-1}R$, and then calculating
$$\tau_* = \tau - \frac{1}{\lambda}.$$\\

\noindent {\tt Step 4.} To check if $\tau_*$ is really optimal, start iterating through the greedy procedure on the ball $\cB_{\tau_*}(A)$, as in {\tt Step 1}.\\

\nid If during some iteration a matrix with negative spectral abscissa is obtained, $\tau_*$ is not the optiaml value. Stop the greedy procedure and return to {\tt Step 1}, taking now $\tau_*$ as the starting value for $\tau$.\\

\nid Else, if we finish the greedy procedure obtaining the matrix $X_*$ with minimal spectral abscissa $\eta(X_*) = 0$ on the ball $\cB_{\tau_*}(A)$, we are done: $\tau_*$ is the optimal value for the problem (\ref{stab.hurwitz}), with $X_*$ as the corresponding optimal solution. {\hfill $\diamond$}\\

\bigskip

\nid \textbf{Remark.} It is worth noting that in Step 1 (and Step 4) of Algorithm 3, we do not need to bring the greedy procedure to its completion. The moment we obtain the matrix with a negative spectral abscissa (not necessarily minimal), we can leave the greedy procedure, thus saving the computational time.\\

\smallskip

\noindent{\bf Example.} For a given Hurwitz unstable Metzler matrix\\

$$
A = 
\left(
\begin{array}{rrrrr}
3  & \phantom{-}0 & \phantom{-}2 & \phantom{-}1 & \phantom{-}4 \\
7  & -4 & 6 & \phantom{-}5  & \phantom{-}7 \\
3  & \phantom{-}4 & \phantom{-}2 & \phantom{-}3 & \phantom{-}0 \\
2  & \phantom{-}1 & \phantom{-}1 & -1 & \phantom{-}8 \\
8  & \phantom{-}0 & \phantom{-}0 & \phantom{-}4 & \phantom{-}9 \\
\end{array} 
\right)
$$\\
\\
with $\eta(A) > 0$, the described procedure computes the matrix\\

$$
X = 
\left(
\begin{array}{rrrrr}
0  & \phantom{-}0 & \phantom{-}0 & \phantom{-}0 & \phantom{-}0 \\
7  & -7 & \phantom{-}6 & \phantom{-}5  & \phantom{-}0 \\
3  & \phantom{-}0 & -4 & \phantom{-}3 & \phantom{-}0 \\
2  & \phantom{-}0 & \phantom{-}0 & -1 & \phantom{-}0 \\
8  & \phantom{-}0 & \phantom{-}0 & \phantom{-}4 & -1 \\
\end{array} 
\right)
$$\\
\\
as the closest stable Metzler with $\tau_* = 10$. 

\bigskip

\begin{center}
\textbf{4.3. Numerical results}
\end{center}

\medskip

In this subsection we report the numerical results for the implementation of Algorithm 3, as the dimension of the starting matrix $A$ is varied. $d$ represents the dimension of the problem, $\#$ the number of iterations required for the procedure to terminate, and $t$ is a computational time.\\

\begin{table}[H]
\begin{center}
\begin{tabular}{c|c c c c c c}
$d$ & 50  & 100 & 250 & 500 & 750 & 1000\\
\hline
$\#$ & 4.9 & 7.6 & 7.3 & 8.3 & 9.9 & 7.4\\
$t$ & 0.19s & 1.08s & 3.06s & 20.78s & 80.98s & 115.03s\\ 

\end{tabular}
\caption*{{\footnotesize Table 4: Results for full Metzler matrices}}
\end{center}
\end{table}

\begin{table}[H]
\begin{center}
\begin{tabular}{c|c c c c c c}
$d$ & 50  & 100 & 250 & 500 & 750 & 1000\\
\hline
$\#$ & 7.1 & 9.6 & 12.6 & 13.4 & 17.1 & 15.4 \\
$t$ & 5.08s & 19.5s & 63.01s & 152.54s & 484.06s & 773.18s\\ 

\end{tabular}
\caption*{{\footnotesize Table 5: Results for sparse Metzler matrices, with 9-15\% sparsity of each row}}
\end{center}
\end{table}

\nid The next table shows how the sparsity affects the computations. The dimension is kept fixed at $d = 850$ in all experiments, while the density parameter $\gamma_i$ for each row of the starting matrix is randomly chosen from the given interval.\\

\begin{table}[H]
\begin{center}
\begin{tabular}{c|c c c c c}
 $\gamma_i$ & 3-8 & 9-15  & 16-21 & 22-51 & 52-76\\
\hline
$\#$ & 15.3 & 15.9 & 10.4 & 21 & 10.2\\
$t$ & 1503.58s & 664.77s & 427.3s & 381.17s & 140.63s\\
\end{tabular}
\caption*{{\footnotesize Table 6: Effects of sparsity}}
\end{center}
\end{table}

\bigskip

\begin{center}
\textbf{4.4. Closest Schur stable non-negative matrix in the $l_\infty$ norm}
\end{center}

\medskip

Building upon the ideas of Algorithm 3, we propose the following modification of the algorithm for Schur stabilization, given in \cite{NePr}.\\

\newpage

\nid {\tt Alg. 4: Computing the closest Schur stable non-negative matrix in $l_\infty$ norm}\\

\bigskip

\nid {\tt Step 0.} Take $\frac{\|A\|}{2}$ as the starting value for $\tau$.\\

\noindent {\tt Step 1.} Start the selective greedy method procedure for minimizing the spectral radius on the ball of non-negative matrices $\cB_\tau^+(A)$. Iterate until a matrix $X$ is obtained with $\rho(X) < 1$. When this is done, stop the greedy procedure, compute its leading eigenvector $\bv$, and proceed to the next step.\\

\nid If greedy procedure finishes finding the Schur unstable matrix as the optimal one on $\cB_\tau^+(A)$, keep implementing the bisection in $\tau$, until the ball $\cB_\tau^+(A)$ containing strongly Schur stable matrix is reached.\\

\noindent {\tt Step 2.} We construct matrices $C = (c_{ij})$ and $R$, as follows. For each $i\in\cS$
\begin{equation*}
c_{ij_k} =
\left\{
\begin{array}{lr}
0, & k < l_i\\
\sum_{s=1}^{l_i} a_{ij_s}, & k = l_i\\
a_{ij_k}, & k > l_i
\end{array}
\right.
\end{equation*}
where $l_i$ the minimal index for which $\sum_{s=1}^{l_i} a_{ij_s} > \tau$. If this $l_i$ does not exist, we put $l_i = m$.  $R$ is a boolean matrix with ones on positions $(i,l_i)$ and zeros in all other places. If some of the indices $i,j$ are not in the support, then we put $c_{ij} = x_{ij}$. We can write $X = C - \tau R$.\\

\nid Denote by $\tau_*$ a potential optimal value of the problem (\ref{stab.schur}). To determine it, proceed to the next step.\\

\noindent {\tt Step 3.} Since $\rho(X) < 1$, we have by Lemma \ref{lem1} that $I - X$ is invertible and
$$(I - X)^{-1} =[I - (C - \tau R)]^{-1} \geqslant 0.$$
Since $1 = \rho(C - \tau_* R) = \rho(C - \tau R + (\tau - \tau_*)R)$, it follows that $\det(I -(C - \tau R) - (\tau - \tau_*)R) = 0$. From here we have\\
\begin{equation*}
\det\Big(\frac{1}{\tau - \tau_*}I - [I - (C-\tau R)]^{-1}R\Big) = 0.
\end{equation*}\\
Matrix $[I - (C - \tau R)]^{-1}R$  is non-negative and $\frac{1}{\tau - \tau_*}$ is its (positive) leading eigenvalue.\\

\nid Therefore, in this step we find the potential optimal value $\tau_*$ by computing the leading eigenvalue $\lambda$ of the matrix $[I -(C - \tau R)]^{-1}R$, and then calculating
$$\tau_* = \tau - \frac{1}{\lambda}.$$\\

\noindent {\tt Step 4.} To check if $\tau_*$ is really optimal, start iterating through the greedy procedure on the ball $\cB_{\tau_*}^+(A)$, as in the {\tt Step 1}.\\

\nid If in some iteration a matrix $Y$ is obtained with $\rho(Y) < 1$, then $\tau_*$ is not optimal. Stop the greedy procedure, return to {\tt Step 1}, and continue doing the bisection taking now $\tau_*$ as the starting value.\\

\nid Else, if we finish the greedy procedure obtaining the matrix $X_*$ with minimal spectral radius $\rho(X_*) = 1$ on the ball $\cB_{\tau_*}^+(A)$, we are done: $\tau_*$ is the optimal value for the (\ref{stab.schur}), with $X_*$ as the corresponding optimal solution. {\hfill $\diamond$}\\

\bigskip

We now present the numerical results for the implementation of Algorithm 4, as the dimension of the starting matrix $A$ is varied.

\begin{table}[H]
\begin{center}
\begin{tabular}{c|c c c c c c}
$d$ & 50  & 100 & 250 & 500 & 750 & 1000\\
\hline
$\#$ & 8.4 & 10.6 & 13.5 & 22.8 & 21.7 & 23.6\\
$t$ & 0.45s & 1.61s & 5.32s & 29.38s & 77.69s & 179.7s\\ 

\end{tabular}
\caption*{{\footnotesize Table 7: Results for totally positive matrices}}
\end{center}
\end{table}

\begin{table}[H]
\begin{center}
\begin{tabular}{c|c c c c c c}
$d$ & 50  & 100 & 250 & 500 & 750 & 1000\\
\hline
$\#$ & 12.7 & 13.1 & 12.5 & 17.5 & 17.1 & 21.7 \\
$t$ & 0.97s & 1.32s & 4.84s & 21.2s & 53.25s & 129.62s\\ 

\end{tabular}
\caption*{{\footnotesize Table 8: Results for sparse non-negative matrices, with 9-15\% sparsity of each row}}
\end{center}
\end{table}

\nid In the next table we shows how the sparsity affects the computations. The dimension is kept fixed $d = 850$ in all experiments, while the density parameter for each row of the starting matrix is randomly chosen from the given interval.\\

\begin{table}[H]
\begin{center}
\begin{tabular}{c|c c c c c}
 $\gamma_i$ & 3-8 & 9-15  & 16-21 & 22-51 & 52-76\\
\hline
$\#$ & 17.8 & 19.7 & 16.5 & 25.2 & 25.1\\
$t$ & 105.22s & 84.25s & 130.35s & 154.02s & 233.16s\\
\end{tabular}
\caption*{{\footnotesize Table 9: Effects of sparsity, for non-negative matrices}}
\end{center}
\end{table}

\bigskip

If we compare these results with the numerical results for Schur stabilization from \cite{NePr}, we can observe a remarkable speed-up in the computational time. The first reason for this is because in the updated algorithm we do not need to conduct the greedy method until the end: we quit it as soon as the strongly Schur stable matrix is obtained. By quitting the greedy method before it finishes, we usually avoid obtaining matrices with zero spectral radius. As practical experiments showed, computation of the leading eigenvector for the zero spectral radius matrices can be drastically slow, especially for the very sparse big matrices. The big chunk of computational time in implementation of the algorithm from \cite{NePr} in fact goes for iterating through zero spectral radius matrices and computing their leading eigenvectors. In the updated procedure this is effectively avoided, providing us with significantly faster computations.\\

\nid Comparing the experimental results for Algorithm 3 and Algorithm 4, we see that Hurwitz stabilization is more time demanding than the Schur stabilization. In Algorithm 3 it is impossible to avoid matrices with spectral abscissa equal to zero, and the computation of their leading eigenvectors can also get very time consuming.          

\bigskip

\begin{center}
\textbf{4.5. The importance of the solution set}
\end{center}

\medskip

The following example illustrates the importance of imposing the set of allowable solutions for the stabilization problems.\\

\nid \textbf{Example.} Observe the non-negative matrix
$$
A = 
\left(
\begin{array}{rr}
1  & 9 \\
6  & 0 \\
\end{array} 
\right)
$$\\
with spectral abscissa (i.e. spectral radius) greater than one. If we want to find a closest \textit{non-negative} matrix having spectral abscissa (i.e. spectral radius) equal to one, we have a Schur stabilization problem. Solving it, we obtain the matrix 
$$
X_* = 
\left(
\begin{array}{ll}
0  & 4.236 \\
0.236  & 0  \\
\end{array} 
\right)
$$
as the closest Schur stable, with $\tau_* = 5.764$. However if we expand our set of admissible solution to a set of Metzler matrices, we get
$$
Y_* = 
\left(
\begin{array}{rr}
-4.4  & \phantom{-}9 \\
\phantom{-}0.6  & \phantom{-}0 \\
\end{array} 
\right)
$$
as the solution. In this case the distance to the starting matrix will be $\tau_* = 5.4$, which is smaller than for the closest non-negative matrix.\\

\nid \textbf{Remark} In the previous example we found the closest Metzler matrix to a matrix $A$ having $\eta(A)=1$ by slightly modifying Step 4: we used greedy procedure on the balls of Metzler matrices $\cB_\tau(A)$, instead on the balls of non-negative matrices $\cB_\tau^+(A)$. Further, we changed spectral radius to spectral abscissa, and determined the indices $l_i$ using the formula (\ref{elle}).\\

\bigskip

\begin{center}\label{sec5}
\large{\textbf{5. Applications}}
\end{center}

\medskip

\begin{center}
\textbf{5.1. Checking the reliability of a given data}
\end{center}

\medskip

One of the most prominent uses of Metzler matrices is for mathematical modelling \cite{FaRi,Lue,Bri}. Numerous phenomena can be modelled by positive linear dynamical systems. In most simple case these systems have the form $\dot{\bx} = A\bx$, where $\bx = \bx(t)$ is a vector of time-dependent unknowns, and $A$ is a time-independent Metzler matrix. The coefficients of $A$ are determined from the gathered data. Often, researchers are interested in examining the stability of the model, which is reflected by the Hurwitz stability of matrix $A$. However, gathered data might contain some errors (i.e. due to the measurement imprecisions), which can lead us to the model with qualitative behaviour completely different than the real picture. To check if our data is reliable, we can optimize spectral abscissa on the ball $\cB_\varepsilon(A)$.\\

For example, let us assume that the matrix of our model-system $A$ is Hurwitz stable and that each entry contains an error not bigger than $\varepsilon$. To check if our model is reliable, we need to maximize the spectral abscissa on $\cB_\varepsilon(A)$.\\

\nid Let $X_*$ be the optimal matrix. If $\eta(X_*) < 0$, we are safe and the qualitative behaviour of our model will agree with the real state of affair, even in the case of making the biggest predicted errors. On the other hand, if we obtain $\eta(X_*) > 0$, we cannot claim that our model reliably describes the phenomenon. In this case one needs to further refine our data gathering methods.\\

Similar reasoning can be applied if for our model we have $\eta(A) > 0$, and we want to check if the real system is unstable as well. However, in this case we would actually need to minimize spectral abscissa on $\cB_\varepsilon$.\\

\bigskip

\begin{center}
\textbf{5.2. Stabilization of 2D postitive linear switching systems}
\end{center}
\medskip

\textit{Positive linear switching systems (LSS)} are an important tool of mathematical modelling. They are a point of an extensive research with many applications \cite{Lib,GuSho,Koz}. Let $A_i,\ i = 1,\ldots,N$ be family of Metzler matrices. A positive LSS is given with:

\begin{equation}\label{LSS}
\begin{array}{lll}
\dot{x}(t) & = & A_\sigma(t)x(t)\\
x(0) & = & x_0,
\end{array}  
\end{equation}

\nid where $x\in\mathbb{R}^d$ and $\sigma: \mathbb{R}_{\geqslant 0}\rightarrow\{1,\ldots,N\}$ is a piecewise constant switching signal.\\

A crucial question in the theory of LSS is the asymptotic stability under the arbitrary switching signal. The following theorem provides us with the necessary condition:

\begin{theorem}{\rm\cite{Lib}}\label{LSS.stable}
If the positive LSS {\rm (\ref{LSS})} is asymptotically stable (under the arbitrary switching), then all the matrices in the convex hull $\textbf{{\rm co}}\{A_1,\ldots,A_N\}$ are Hurwitz stable.
{\hfill $\Box$}
\end{theorem}

\nid The converse of Theorem \ref{LSS.stable} is true only for the two-dimensional case~\cite{GuSho}. Therefore, we can use Hurwitz stabilization on the unstable convex combinations to build an algorithm for the stabilization of the 2D positive LSS. By this we mean constructing a stable LSS from the original system, as shown below.\\

Assume that $2\times2$ Metzler matrices $A_i,\ i = 1,\ldots,N$ are Hurwitz stable. If some of them is not, find and replace it with its closest stable. Suppose that there exists a matrix in $\textbf{{\rm co}}\{A_1,\ldots,A_N\}$ that is not Hurwitz stable.\\

\nid Let $A=\sum_{i=1}^N \alpha_i A_i$ be a convex combination in $\textbf{{\rm co}}\{A_1,\ldots,A_N\}$ with the largest spectral abscissa $(\sum_{i=1}^N \alpha_ 1 = 1,\ \alpha_i\in[0,1],\ i=1,\ldots,N)$. Find its closest Hurwitz stable matrix $A'$, and denote by $\tau_*$ the optimal distance.\\

\nid We now need to decompose the matrix $A'=\sum_{i=1}^N \alpha_i A'_i$, while taking care about matrices $A_i$. We have
\begin{equation}\label{inequa}
\tau_* = \|A - A'\|\leqslant\sum_{i=1}^N \alpha_i\|A_i - A'_i\|=\sum_{i=1}^N \alpha_i\tau_i.
\end{equation}
\nid Choose values $\tau_i$ and matrices $A'_i\in\cB_{\tau_i}(A_i)$, so that (\ref{inequa}) and $A' = \sum_{i=1}^N \alpha_i fA'_i$ are satisfied.\\

\nid Let $A_*$ be a convex combination in $\textbf{{\rm co}}\{A'_1,\ldots,A'_N\}$. If $\eta(A_*) < 0$, then we are done: the switching system built from matrices $A'_i$ is asymptotically stable. Else, we should make a different choice of $\tau_i$, or restart everything, but now starting with matrices $A'_i$.\\

\newpage

\begin{center}
\textbf{5.3. Closest stable sign-matrix and its application to positive LSS}
\end{center}
\medskip

The notions of sign-matrices and sign-stability originated from the problems in economy \cite{QuRu, Sam}, ecology \cite{May, Lev}, and chemistry \cite{Cla}. Ever since, those concepts have been a point of interest in mathematical literature \cite{Bri, JeKl, MayQu}. Metzler sign-matrices are very useful for the dynamical system modelling. They come in handy if we do not posses quantitative data, but just the information on the sign of coefficients of the observed system. By analysing the sign-matrix of the system, we can discern its qualitative behaviour. As we will see later on, Metzler sign-matrices can be used as a tool for stability analysis of linear switching systems.\\

Denote by $\cM_\sgn$ set of all real Metzler matrices with entries from  $\{-1,0,1\}$. Hurwitz stable matrices from $\cM_\sgn$ have one peculiar property: replacing any entry by an arbitrary real number of the same sign does not influence the stability.\\

\nid \textbf{Example.} The matrix $A\in\cM_\sgn$\\

$$
A = 
\left(
\begin{array}{rrrrr}
-1 & 1 & 0 & 0 & 0\\
0 & -1 & 0 & \phantom{-}0 & \phantom{-}1\\
0 & \phantom{-}1 & -1 & \phantom{-}0 & \phantom{-}0 \\
1 & \phantom{-}1 & \phantom{-}0 & -1 & \phantom{-}1 \\
0 & \phantom{-}0 & \phantom{-}0 & \phantom{-}0 & -1 \\
\end{array} 
\right)
$$\\

\nid is Hurwitz stable with $\eta(A) = -1$. Metzler matrix $A'$ with the same sign pattern, given by\\

$$
A' = 
\left(
\begin{array}{ccccc}
 -10^{-9} & 10^9 & \phantom{-}0\phantom{^{-9}} & \phantom{-}0\phantom{^{-9}} & \phantom{-}0\phantom{^{-9}}\\
\phantom{-}0\phantom{^{-9}} & -10^{-9} {\ } & \phantom{-}0\phantom{^{-9}} & \phantom{-}0\phantom{^{-9}}  & 10^9\\
\phantom{-}0\phantom{^{-9}} & 10^9 & -10^{-9} & \phantom{-}0\phantom{^{-9}} & \phantom{-}0\phantom{^{-9}}\\
10^9 & 10^9 & \phantom{-}0\phantom{^{-9}} & -10^{-9} & 10^9\\
\phantom{-}0\phantom{^{-9}} & \phantom{-}0\phantom{^{-9}} & \phantom{-}0\phantom{^{-9}} & \phantom{-}0\phantom{^{-9}} & -10^{-9} {\ }\\
\end{array} 
\right)
$$\\

\nid has $\eta(A') = -10^{-9}$, i.e., is also Hurwitz stable. It remains stable in spite of very 
big changes ($\cO(10^9)$) of all non-negative entries!\\
\nid This property can be derived by studying sign-matrices, which we introduce now.

\begin{defi}
A matrix is called a {\rm sign-matrix} if its entries take values from the set $\{-,0,+\}$. If $-$ appears only on the main diagonal, we say it is {\rm Matzler sign-matrix}. 
\end{defi}

\nid 

\nid In analogy with the real Metzler matrices, we can also consider Hurwitz stability of sign-matrices.

\begin{defi}
Let $M$ be a sign-matrix. We say that real matrix $X$ belongs to a {\rm qualitative class of a matrix} $M$, denoted by $\cQ(M)$, if
\begin{equation*}
x_{ij}
\left\{
\begin{array}{ll}
< 0, & m_{ij} = -\\
= 0, & m_{ij}= 0\\
> 0, & m_{ij} = +.
\end{array}
\right.
\end{equation*}
\end{defi}

\begin{defi}
A Metzler sign-matrix $M$ is {\rm (strongly) Hurwitz stable} if all the matrices from the qualitative class $\cQ(M)$ are (strongly) Hurwitz stable.
\end{defi}

\begin{theorem}{\rm \cite{Bri}}\label{thm10}
Let $M$ be a Metzler sign-matrix, and $\sgn(M)$ be a real matrix given by
\begin{equation*}
m^\sgn_{ij} = 
\left\{
\begin{array}{ll}
-1, & m_{ij} = -,\\ 
0, & m_{ij} = 0,\\
1, & m_{ij} = +.
\end{array}
\right.
\end{equation*}
\nid $M$ is (strongly) Hurwitz stable if and only if $\sgn(M)$ is (strongly) Hurwitz stable. Moreover, $M$ is strongly Hurwitz stable if and only if $\sgn(M) + I$ is the adjacency matrix of an acyclic graph.
{\hfill $\Box$}
\end{theorem}

\nid From Theorem \ref{thm10} we see why the Hurwitz stable matrices from $\cM_\sgn$ do not change their stability, even as we change their non-zero entries.\\

The problem of Hurwitz stabilization of real Metzler matrices can be formulated for the Metzler sign-matrices as well. Theorem \ref{thm10} gives a way to do so: for a given Hurwitz unstable Metzler sign-matrix $M$, we need to find \textit{closest Hurwitz stable Metzler sign-matrix} $X$. In other words, $X$ should solve
\begin{equation}\label{closest.hurwitz.sign}
\left\{
\begin{array}{l}
\|\sgn(M) - \sgn(X)\|\ \rightarrow\min\\
\eta(\sgn(X))\leqslant 0.
\end{array}
\right.
\end{equation}

\medskip

\nid \textbf{Remark.} Notice that in Hurwitz stabilization of sign-matrix we allow our optimal solution to have negative spectral abscissa as well. We do this since in some cases a solution with zero spectral abscissa does not exist\footnote{\ see Examples below.}. \\ 

Denote $\cB_k^\sgn =  \cB_k\cap\cM_\sgn$. Using minimization of spectral abscissa on $\cB_k^\sgn$ we can present a simple procedure for solving problem (\ref{closest.hurwitz.sign}).\\

\nid {\tt Algorithm 5: Finding the closest Hurwitz stable sign-matrix}\\

\nid Let $M$ be a Hurwitz unstable sign-matrix. Set $k = \lfloor{\frac{\|\sgn(A)\|}{2}}\rfloor$. By $X_k$ denote the sign-matrix such that $\eta_k = \eta(\sgn(X_k))$ is minimal on $\cB_k^\sgn$.\\

\nid Doing a bisection in $k$, minimize the spectral abscissa on the ball $\cB_k^\sgn$. Do this until sign-matrices $X_{k-1}$ and $X_k$ are obtained, with the (minimal) spectral abscissas $\eta_{k-1} > 0$ and $\eta_{k}\leqslant 0$.\\

\nid Take $k = k_*$ as the optimal value and $X_k = X_*$ as the optimal solution of (\ref{closest.hurwitz.sign}). {\hfill $\diamond$}\\

\medskip

\nid \textbf{Example.} For the unstable Metzler sign-matrix 

$$
M = 
\left(
\begin{array}{ccccc}
 0 &  + &  + & + & 0\\
 + &  + &  0 & + & +\\
 + &  + &  0 & 0 & +\\
 + &  0 &  0 & - & +\\
 0 &  0 &  + & + & +\\
\end{array} 
\right)
$$

\nid we find 

$$
M_* = 
\left(
\begin{array}{ccccc}
 0  & 0 & 0 & + & 0\\
 +  & - & 0 & + & +\\
 +  & 0 & - & 0 & +\\
 0  & 0 & 0 & - & 0\\
 0  & 0 & 0 & + & 0\\
\end{array} 
\right)
$$\\

\nid with $\eta(\sgn(M_*)) = 0$ as the closest sign-stable, and optimal distance $k_* = 2$.\\

\smallskip

\nid \textbf{Example.} For the unstable Metzler sign-matrix 

$$
M = 
\left(
\begin{array}{ccccc}
 - & + & 0 & 0 & +\\
 + & 0 & 0 & + & +\\
 + & 0 & 0 & + & 0\\
 + & + & + & 0 & 0\\
 0 & + & + & 0 & -\\
\end{array} 
\right)
$$

\nid we find

$$
M_* = 
\left(
\begin{array}{ccccc}
 - & 0 & 0 & 0 & 0\\
 + & - & 0 & 0 & +\\
 + & 0 & - & 0 & 0\\
 + & 0 & + & - & 0\\
 0 & 0 & 0 & 0 & -\\
\end{array} 
\right)
$$\\

\nid with $\eta(\sgn(M_*)) = -1$ as the closest sign-stable, and optimal distance $k_* = 2$. The optimal solution with zero spectral abscissa is impossible to find, because at the distance $k = 1$ we obtain the matrix

$$
M' = 
\left(
\begin{array}{ccccc}
 - & 0 & 0 & 0 & +\\
 + & 0 & 0 & 0 & +\\
 + & 0 & 0 & 0 & 0\\
 + & + & + & - & 0\\
 0 & 0 & + & 0 & -\\
\end{array} 
\right)
$$\\

\nid with $\eta(\sgn(M')) = 0.46$ as the one with the minimal spectral abscissa.\\ 

\medskip

We have seen that Theorem \ref{LSS.stable} works both ways only in the two-dimensional case. However, there is a criterion of  asymptotic stability, valid in any dimension, which involves sign-matrices and sign-stability.

\begin{theorem}{\rm \cite{Bri}}\label{switch.sign}
Let $M_i,\ i=1,\ldots,N$ be a family of Metzler sign-matrices. For all $A_i\in\cQ(M_i),\ i=1,\ldots,N$,  the LSS {\rm (\ref{LSS})} is asymptotically stable if and only if all the diagonal entries of $M_i,\ i=1,\ldots,N$ are negative and $\sum_{i=1}^N M_i$ is Hurwitz sign-stable.
{\hfill $\Box$}
\end{theorem}

\nid Now, assume that the sign-matrices $M_i,\ i=1,\ldots,N$ have negative diagonal entries, that they are Hurwitz stable, but their sum $M = \sum_{i=1}^N M_i$ is not. Using Algorithm 5, we can find closest Hurwitz sign-stable matrix $M'$ to $M$. Let $k_*$ be the optimal value. Now we need to decompose $M'$ into the sum $\sum_{i=1}^N M'_i$, while taking into the consideration the structure of matrices $M_i$. We have
\begin{equation}\label{taustar}
k_* = \|M - M'\| \leqslant \sum_{i=1}^N \|M_i - M'_i\| = \sum_{i=1}^N k_i.
\end{equation}
\nid We need to choose values $k_i$ and matrices $M'_i\in\cB_{k_i}^\sgn(M_i)$, so that (\ref{taustar}) and $M' = \sum_{i=1}^N M'_i$ are satisfied. Matrices $M'_i$, including their sum, will have negative entries on the main diagonal and be Hurwitz stable. Therefore, by Theorem \ref{switch.sign}, any LSS we construct from matrices in $\cQ(M'_i)$ will be asymptotically stable, while keeping the similar structure to the initial system.\\

\medskip

\nid\textbf{Example.} Observe a positive linear switching system, switching between the following three strongly Hurwitz stable matrices:

$$
A_1 = 
\left(
\begin{array}{rrrr}
-8 & \phantom{-}0 & \phantom{-}0 & \phantom{-}0\\
\phantom{-}0 & -4 & \phantom{-}2 & \phantom{-}0\\
\phantom{-}7 & \phantom{-}0 & -8 & \phantom{-}3\\
\phantom{-}0 & \phantom{-}1 & \phantom{-}9 & -9\\
\end{array} 
\right)\ 
A_2 = 
\left(
\begin{array}{rrrr}
-5 & \phantom{-}0 & \phantom{-}0 & \phantom{-}2\\
\phantom{-}0 & -3 & \phantom{-}0 & \phantom{-}0\\
\phantom{-}2 & \phantom{-}2 & -1 & \phantom{-}3\\
\phantom{-}0 & \phantom{-}5 & \phantom{-}0 & -8\\
\end{array} 
\right)\
A_3 = 
\left(
\begin{array}{rrrr}
-9 & \phantom{-}0 & \phantom{-}2 & \phantom{-}0\\
\phantom{-}0 & -9 & \phantom{-}0 & \phantom{-}2\\
\phantom{-}0 & \phantom{-}7 & -2 & \phantom{-}3\\
\phantom{-}0 & \phantom{-}0 & \phantom{-}0 & -4\\
\end{array} 
\right).
$$\\

\nid Their corresponding sing-matrices are given by

$$
M_1 = 
\left(
\begin{array}{cccc}
- & 0 & 0 & 0\\
0 & - & + & 0\\
+ & 0 & - & +\\
0 & + & + & -\\
\end{array} 
\right)\
M_2 = 
\left(
\begin{array}{cccc}
- & 0 & 0 & +\\
0 & - & 0 & 0\\
+ & + & - & +\\
0 & + & 0 & -\\
\end{array} 
\right)\
M_3 = 
\left(
\begin{array}{cccc}
- & 0 & + & 0\\
0 & - & 0 & +\\
0 & + & - & +\\
0 & 0 & 0 & -\\
\end{array} 
\right),
$$\\

\nid and their sum $M = M_1 + M_2 + M_3$ with

$$
M = 
\left(
\begin{array}{cccc}
- & 0 & + & +\\
0 & - & + & +\\
+ & + & - & +\\
0 & + & + & -\\
\end{array} 
\right).
$$\\

\nid Since sign-matrix $M$ is Hurwitz unstable (with $\eta(\sgn(M)) = 1.303$), our positive LSS will not be asymptotically stable. Applying Algorithm 6 on it, we find its closest strongly Hurwitz stable sing-matrix

$$
M' = 
\left(
\begin{array}{cccc}
- & 0 & 0 & +\\
0 & - & 0 & +\\
+ & + & - & 0\\
0 & + & 0 & -\\
\end{array} 
\right),
$$\\

\nid with an optimal distance $k_* = 1$. We select $k_1 = k_2 = k_3 = 1$ and decompose $M' = M'_1 + M'_2 + M'_3$ by choosing $M'_i\in\cB_{k_i}^\sgn(M_i),\ i=1,2,3$. We have:

$$
M'_1 = 
\left(
\begin{array}{cccc}
- & 0 & 0 & 0\\
0 & - & 0 & 0\\
+ & 0 & - & +\\
0 & + & 0 & -\\
\end{array} 
\right)\
M'_2 = 
\left(
\begin{array}{cccc}
- & 0 & 0 & +\\
0 & - & 0 & 0\\
+ & + & - & 0\\
0 & + & 0 & -\\
\end{array} 
\right)\
M'_3 = 
\left(
\begin{array}{cccc}
- & 0 & 0 & 0\\
0 & - & 0 & +\\
0 & + & - & 0\\
0 & 0 & 0 & -\\
\end{array} 
\right).
$$\\

\nid Following the sign pattern of matrices $M'_i$ and correspondingly cutting the interdependencies in matrices $A_i$, we get

$$
A'_1 = 
\left(
\begin{array}{rrrr}
-8 & \phantom{-}0 & \phantom{-}0 & \phantom{-}0\\
\phantom{-}0 & -4 & \phantom{-}0 & \phantom{-}0\\
\phantom{-}7 & \phantom{-}0 & -8 & \phantom{-}0\\
\phantom{-}0 & \phantom{-}1 & \phantom{-}0 & -9\\
\end{array} 
\right)\ 
A'_2 = 
\left(
\begin{array}{rrrr}
-5 & \phantom{-}0 & \phantom{-}0 & \phantom{-}2\\
\phantom{-}0 & -3 & \phantom{-}0 & \phantom{-}0\\
\phantom{-}2 & \phantom{-}2 & -1 & \phantom{-}0\\
\phantom{-}0 & \phantom{-}5 & \phantom{-}0 & -8\\
\end{array} 
\right)\
A'_3 = 
\left(
\begin{array}{rrrr}
-9 & \phantom{-}0 & \phantom{-}0 & \phantom{-}0\\
\phantom{-}0 & -9 & \phantom{-}0 & \phantom{-}2\\
\phantom{-}0 & \phantom{-}7 & -2 & \phantom{-}0\\
\phantom{-}0 & \phantom{-}0 & \phantom{-}0 & -4\\
\end{array} 
\right).
$$\\

\nid A newly obtained LSS built from matrices $A'_i$ will be asymptotically stable under the arbitrary switching.\\

\nid\textbf{Acknowledgement.} Countless thanks to prof V. Yu. Protasov for the hours of inspiring discussions, all the valuable remarks and the support he was selflessly giving during the writing of this paper.

\end{document}